\documentclass[12pt]{article} 
\usepackage{graphicx}
\usepackage{amsthm}
\usepackage{amssymb}
\usepackage{fdsymbol}
\usepackage{amsmath}
\usepackage{mathrsfs}
\newtheorem{theorem}{Theorem}[section]
\newtheorem{proposition}{Proposition}[section]
\newtheorem{corollary}{Corollary}[section]
\newtheorem{lemma}{Lemma}[section]
\newtheorem{definition}{Definition}[section]

\newtheorem{remark}{Remark}[section]
\newtheorem{notation}{Notation}[section]
\newtheorem{example}{Example}[]

\begin{document}

\title{Periodicity and Free Periodicity of Alternating Knots}
\author{Antonio F. Costa and Cam Van Quach Hongler }

\maketitle

\begin{abstract} 
A knot $K$ in $S^3$ is {\bf $q$-periodic} if it admits a symmetry that is conjugate to a rotation of order $q$ of $S^3$. If $K$ admits a symmetry which is a homeomorphism without fixed points of period $q$ of $S^3$, then $K$ is called {\bf freely $q$-periodic}.

In a previous paper \cite{co2}, we obtained, as a consequence of  Flyping Theorem due to Menasco and Thislethwaite, that the $q$-periodicity with $q>2$ can be visualized in an alternating projection as a rotation of the projection sphere. See also \cite{boy}.

In this paper, we show that the free $q$-action of an alternating knot can be represented on some alternating projection as a composition of a rotation of order $q$ with some flypes all occurring on the same twisted band diagram of its essential Conway decomposition. Therefore, for an alternating knot to be freely periodic, its essential decomposition must satisfy certain conditions. We show that any free or non-free $q$-action is somehow visible ({\bf virtually visible}) and give some sufficient criteria to detect the existence of $q$-actions from virtually visible projections.\\
Finally, we show how the Murasugi decomposition into atoms as initiated in \cite{quwe1} and \cite{stoi} enables us to determine the visibility type $(q,r)$ of the $q$-periodic freely alternating knots ($(q,r)$-lens knots \cite{chbi}); in fact, we only need to focus on a certain atom of their Murasugi decomposition in order to deduce their visibility type.
\end{abstract}

\section{Introduction}
In this paper, links (knots are one-component links) in $S^3$ and projections in $S^2$ are assumed to be oriented and indecomposable. By a {\bf projection} on $S^2$, we mean the image of a link in $S^3$ by a generic projection onto the 2-sphere,  denoted $\mathbb S$, where crossings are labeled to reflect under and over arcs at crossings. \\
The visibility of the $q$-periodicity of oriented alternating links is studied in \cite{co}.\\
We now recall the notion of visibility of the $q$-periodicity of an alternating knot as defined in \cite{co2}.
 \begin{definition} Let $K$ be an alternating $q$-periodic knot. The $q$-periodicity of $K$ is {\bf visible} if $K$ has an alternating projection that displays its $q$-periodicity as a $2\pi \over q$-rotation. Such an alternating projection is  called a {\bf $q$-periodic projection}.
 \end{definition}
 
 In \cite{co2}, it is shown that:

\begin{theorem}(Visibility Theorem for $q$-periodicity) Let $K$ be a prime oriented alternating knot that is $q$-periodic with $q \geq 3$. Then $K$ has a $q$-periodic alternating projection.
\end{theorem}

In fact, Theorem 1.1 can be considered as a corollary of the following theorem:

 {\bf Theorem 5.3} {\it (Visibility Theorem for $q$-periodicity with $q \geq 3$).
 Let $K$ be a prime oriented alternating knot which is $q$-periodic with $q \geq 3$. Then there exists an alternating projection $\Pi$ of $K$ such that the $q$-periodic homeomorphism of pairs $\Phi:(S^3,K) \rightarrow (S^3,K)$ is topologically conjugate to an isomorphism $ \Phi_\Pi$ of a realized projection $(S^3, \lambda (\Pi))$ onto itself which is reduced to a flat homeomorphism $\phi$ of order $q$ 
  $$\Phi_\Pi =\phi$$
The restriction $\phi ^*$ of $\phi$ on the projection sphere $\mathbb S$ is a rotation of angle $2 \pi \over q$ without fixed points on the projection $\Pi$.} 

The definitions of a {\it realized projection} and of an {\it isomorphism of realized projections} are given in \S 5.1.

According to the definition of the $q$-periodic visibility for an alternating knot given in \cite{co2}, L. Paoluzzi (\cite{pao}) noticed that the free $q$-periodicity with $q \geq 3$ cannot be visible since the free $q$-periodic map (Definition 2.3) cannot leave invariant a projection 2-sphere. The natural question is how to express the nearest notion of visibility for the free $q$-periodicity in the class of alternating knots in light of the Menasco-Thislethwaite Flyping Theorem:
 
\begin{theorem} (Flyping Theorem \cite {meth}) Let $\Pi_1$ and $\Pi_2$ be two reduced alternating projections of links. If $f:(S^3, \lambda(\Pi_1)) \rightarrow (S^3, \lambda(\Pi_2))$ is a homeomorphism of pairs, then $f$ is a composition of non-trivial flypes and flat homeomorphisms.
\end{theorem}
The aim of this paper is to extend in a way the Visibility Theorem 1.1 of $q$-periodicity to the free $q$-periodicity case in the class of alternating knots. 

By the Geometrization Theorem due to W. Thurston and G. Perelman, it is known that a $q$-freely periodic link in $S^3$ is topologically conjugate to a $(q,r)$-lens link  (Definition 4.3) for some integer $r$ where $0 < r < q$ and $gcd (r,q)=1$. 

A $(q,r)$-lens link in $S^3$ is a link that is invariant under the  $(q,r)$-action $\Phi_{q,r}$ of $ \mathbb Z_q$ on $S^3=\{(z_1,z_2) \in \mathbb C \times \mathbb C \, | \, |z_1|^2 +|z_2|^2=1\}$ such that $\Phi_{q,r}$ is defined by the following Equation (1):
\begin{align*}
\Phi_{q,r}&:S^3 \quad \longrightarrow  S^3 \\
&(z_1,z_2)  \longmapsto (e^{i2\pi\over q}z_1,e^{i2\pi r\over q}z_2)
\end{align*}

One of the principal tools of our study is the decomposition of an alternating projection by its set of {\it essential} (resp. {\it canonical}) {\ it Conway circles} into {\it jewels} and {\it twisted band diagrams} which leads to the construction of its {\it essential }(resp. {\it canonical}) {\it structure tree} ({\cite{co2}). It is a consequence of Flyping Theorem that the essential structure tree such as the canonical structure tree is an invariant in the class of alternating links.\\
In the case  $r=0$, $\Phi_{q,0}$ induces a $q$-action whose set of fixed points is the circle $l_2=\{(0,z_2) \in S^3 \}$. The links that do not cut  $l_2$ and that are invariant under the $q$-action are called $q$-periodic.

The free $q$-periodicity can be studied in a similar way as the $q$-periodicity which has been treated in \cite{co2} using the essential decomposition for alternating projections. Let $K$ be a prime alternating knot in $S^3$ with its essential structure tree $\widetilde{\cal{A}} (K)$. Let us assume that $K$ has a $q$-action which induces the automorphism $\widetilde \varPhi$ on $\widetilde{\cal{A}} (K)$. In the case where the fixed point of $\widetilde \varPhi$ is a jewel, $K$ is necessarily $q$-periodic and not freely $q$-periodic. The relevant case is when $Fix(\widetilde \varPhi)$ is a vertex $V_0$ corresponding to a twisted band diagram called the {\bf main twisted band diagram} which is invariant by the isomorphism of the realized projections  $\Phi_\Pi$ (defined in \S5.1).
 
For the case of free $q$-periodicity, we have: 

{\bf Theorem 5.4} {\it (Visibility Theorem for free $q$-periodicity with $q \geq 3$) 
 Let $K$ be a prime oriented alternating knot with an action defined by a freely $q$-periodic homeomorphism of pairs $\Phi:(S^3,K) \rightarrow (S^3,K)$ where $q \geq 3$.
Then there exists an alternating projection $\Pi$ of $K$ with a twisted band diagram $\mathcal T_0$ such that $\Phi$ is topologically conjugate to an isomorphism of realized projections  $\Phi_\Pi: (S^3, \lambda (\Pi)) \rightarrow (S^3, \lambda (\Pi))$ which satisfies the following conditions:
$$
 \Phi_\Pi =\phi \circ F_0
$$
where \\
(1) $\phi$ is a flat homeomorphism of order $q$ whose principal part is topologically conjugate to a rotation, leaving invariant the diagram $\mathcal T_0$ and with generic orbits in its action on the other diagrams of the essential Conway decomposition, and\\
(2) $F_0$ is a non-trivial composition of standard flypes on $\lambda (\Pi)$ all occurring on $\lambda(\mathcal T_0)$.}

We call $\mathcal T_0$ the {\bf main twisted band diagram} of $\Pi$.
 
Let $K$ be a freely periodic alternating knot. For any alternating projection $\Pi$ of $K$, we cannot reduce the composition of flypes $F$ in $ \Phi_\Pi$ of Flyping Theorem 5.1 to the identity map. We will show in \S 5.5 that the best possible visibility we can get for a free $q$-action is obtained as described in Theorem 5.4 quoted above.

We introduce the notion of {\it virtual $q$-visibility} (see Definition 5.6) for an alternating projection which concerns both cases of $q$-periodicity and free $q$-periodicity. From a $q$-virtually visible projection of a freely $q$-periodic alternating knot $K$, we can deduce the homeomorphism $\Phi_{q,r}$ of $S^3$ which characterizes $K$ as a $(q,r)$-lens link. Using virtually visible projections, we also present some criteria for detecting free and non-free $q$-periodicity.
 
In \S 6, by applying the Murasugi decomposition into {\it atoms} as initiated by \cite{quwe1} and \cite{stoi} for oriented prime alternating knots $K$ endowed with a $q$-action, we can describe the automorphism induced by its $q$-action on its {\it adjacency graph} ${\cal G}(K)$ (\cite{co}). It is interesting to realize that the free q-periodic map $\Phi_{q,r}$ is completely characterized by the main Murasugi atom corresponding to the vertex of ${\cal G}(K)$ invariant under this automorphism.
To derive the pair $(q,r)$ which characterizes the free $q$-map of $K$, it is sufficient to focus on the {\it main Murasugi atom} corresponding to the invariant vertex under the automorphism induced on ${\cal G}(K)$.

We would like to thank Professor Mar\'{\i}a Teresa Lozano for some comments and corrections to the first version of this article. The first author is partially supported by PGC2018-096454-B-100 (Spanish Ministry of Science, Innovation and Universities).
 
 \section{Free and non-free $q$-actions}
We first recall the definitions of $q$-periodicity and free $q$-periodicity of a knot where $q>2$.

\begin{definition}
Let $q$ be an integer $>2$. A  knot $K$ is $q$-periodic if there is a homeomorphism of pairs $\Phi$ from $(S^3, K)$ to $(S^3, K)$ of period $q$ that satisfies the following conditions:\\
(1) $\Phi$ is topologically conjugate to a $2\pi \over q$-rotation around a ``line" (circle) $\alpha$ in $S^3$. \\
(2) $\alpha \cap K= \emptyset$.\\
$\Phi$ defines a {\bf non-free $q$-action} on $(S^3, K)$ and $\Phi$ is called a {\bf $q$-periodic map} of $K$.
\end{definition}
 
\begin{definition} 
A knot is strictly $q$-periodic if it is $q$-periodic but not $q \over s$-periodic for any $s \geq 2$.
\end{definition}
  
\begin{definition} Let $q$ be an integer $ > 2$. A knot $K$ in $S^3$ is freely $q$-periodic if there is a homeomorphism of pairs $\Phi$ from $(S^3,K)$ to $(S^3,K)$ of period $q$ such that $\Phi^i$ has no fixed point for all $1\leq i \leq q-1$.\\
$\Phi$ defines a {\bf free $q$-action} on  $(S^3,K)$ and $\Phi$ is called a {\bf free $q$-periodic map} of $K$.  
\end{definition}

From now on, by a (free) $q$-periodic knot (map), we mean a strictly $q$-periodic knot (map).

\section{Canonical Decomposition of a Projection}
In this section, we briefly recall the notions on the {\bf canonical decomposition of a projection} of a link (\cite{quwe0}). For more details and in particular on the essential structure trees deduced from the canonical decomposition of alternating projections, we refer to \cite{co2}.

{\bf For the canonical decomposition, the orientation of the projection does not matter}.

A {\bf diagram} of a link projection $\Pi$ is a pair $D = (\Sigma , \Gamma= \Pi \cap \Sigma)$ where $\Sigma$ is a compact connected planar surface embedded on the projection sphere $\mathbb{S}$. 

A Haseman circle of a diagram $D = (\Sigma, \Gamma)$ is a circle $\gamma \subset \Sigma$ that intersects the projection $\Pi$ exactly in 4 points away from the crossings.

\begin{definition}
A {\bf twisted band diagram} (abbreviated as TBD) is a diagram $(\Sigma, \Gamma)$ that is homeomorphic to the diagram of Fig.1 where $\gamma_1 , \dots , \gamma_{k+1}$ with $k \geq 0$ denoting the $k+1$ boundary components of $\Sigma$  are Haseman circles. With the definition of the {\bf weight of a crossing}  given in Fig.2, $a_i$ denotes the total weight of crossings between $\gamma_i$ and $\gamma_{i+1}$.
  
  A {\bf spire} is a TBD with $k+1=1$ boundary component. A {\bf twisted annulus} is a TBD with $k+1=2$ components.

\end{definition}
\begin{figure}[ht]    
  \centering
\includegraphics[scale=.4]{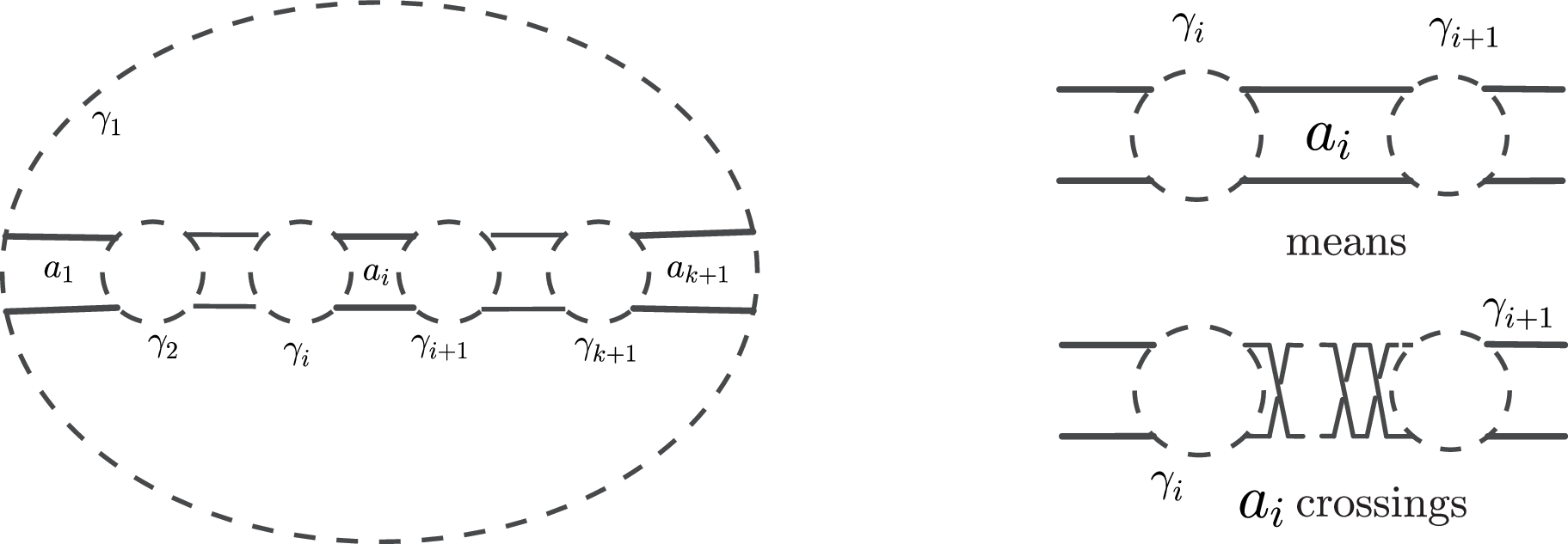}
\caption{ A twisted band diagram}
\end{figure}

\begin{figure}[ht]    
  \centering
   \includegraphics[scale=0.3]{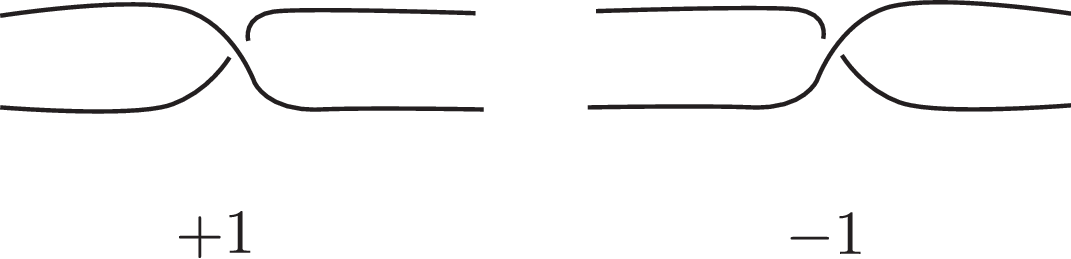}
\caption{The weight of a crossing on a band}
\end{figure}

\begin{definition}
 A {\bf jewel} is a diagram $J$ such that:\\
(1) $J$ is not a twisted band diagram with $k+1=2$ and $a=\pm1$ or with $k+1=3$ and $a=0$ and\\
(2) each Haseman circle of $J$ is boundary-parallel (see \cite{co2}).
 \end{definition}

A {\bf twist region} of a TBD $T_K$  is a portion of $T_K$  between two consecutive circles $\gamma_i$ and $\gamma_{i+1}$ that are  boundaries of  $T_K$ .\\
We now recall the notion of a {\bf flype} which is an essential transformation in our study.
\begin{definition} A flype is a transformation of the projection of a link as described by Fig.3. The crossing in Fig.3 which is moved by the flype is an {\bf active crossing}.
  \end{definition}
 \begin{figure}[h!]    
  \centering
   \includegraphics[scale=.3]{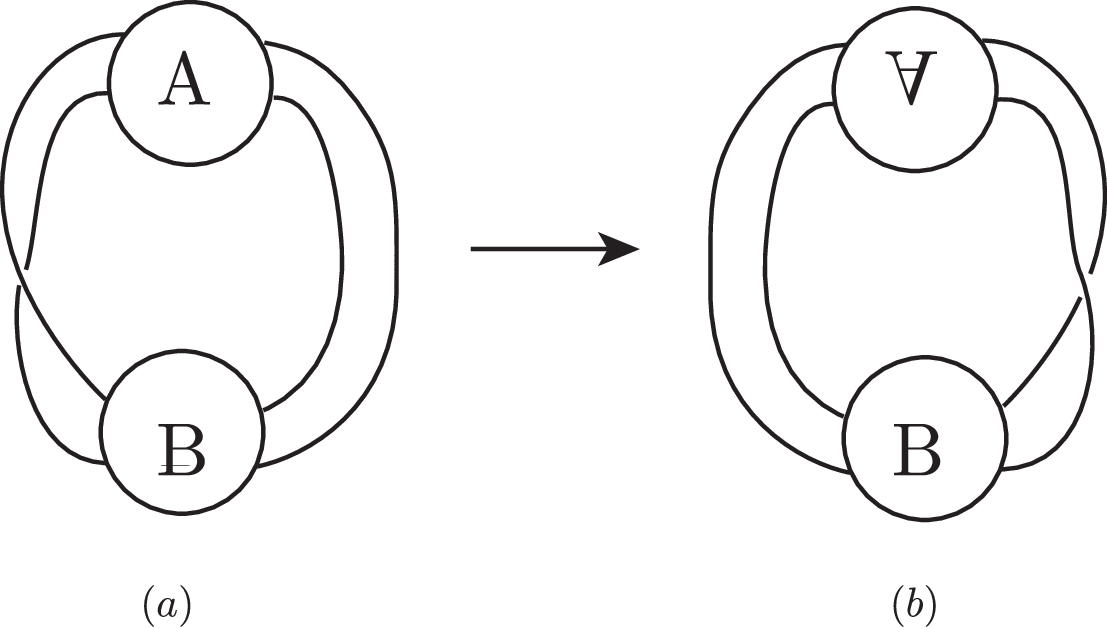}
\caption{A flype}
\end{figure}
By using flypes and Reidemeister moves of type II, we can reduce the number of crossings of a TBD so that all its crossings have  the same sign. Thus, we can assume that all $a_i$ shown in Fig.1 have the same sign and we can define the {\bf weight of a TBD} as  $w=\sum_{i=1}^{k+1}a_i$. 

\subsection{Canonical Conway and Essential Conway circles}

If not otherwise specified, the projections we consider are connected, indecomposable and alternating.

Before defining the canonical and essential Conway circles, we recall the notion of a tangle.

\begin{definition}
A two-dimensional {\bf 2}-tangle $\mathcal{T}$ of a projection $\Pi $ is a diagram $\mathcal{T}=(\Delta , \tau _{\Delta })$ where $\Delta $ is a disk in the projection sphere $\mathbb{S}$ and where the 1-manifold $\tau _{\Delta }=\Pi \cap \Delta$ intersects the boundary $\partial \Delta $ of the disk $\Delta $ at 4 points SE, NE, NW, and SW; these points are located southeast, northeast, and so on $\partial \Delta$.\\
The {\bf boundary} of $\cal T$ is denoted by ${\partial \cal{T}}=\partial \Delta $.
\end{definition} 
Unless otherwise specified, we abbreviate ``two-dimensional tangle" by ``tangle" and a
{\bf 2}-tangle ({\bf 2} means here ``with 2 strands").\\

\begin{definition}
Let $\mathcal{T}=(\Delta ,\tau _{\Delta })$ and $\mathcal{T}{^{\prime }}=(\Delta ,\tau_{\Delta }^{\prime })$ be two 2-tangles of $\Pi$ with the same four endpoints SE, NE, NW and SW. They are \textbf{isotopic} if we can change $\mathcal{T}$ into $\mathcal{T}{^{\prime }}$ by a sequence of Reidemeister moves inside $\Delta$ while keeping the four endpoints fixed.
\end{definition} 

\subsubsection {Canonical Conway circles}
 There exists, up to isotopy respecting $\Pi$, a unique {\it minimal admissible} family of Haseman circles denoted ${\cal C}_{can}(\Pi)$ which decompose $\Pi$ into TBDs and jewels (see Theorem 1 in \cite{quwe0}). We call it the family of {\bf canonical Conway circles} of $\Pi$. 
Recall that \\
An element of the canonical Conway family can be of 3 types:
\begin{enumerate}
\item a circle that separates two jewels.
\item a circle that separates two twisted band diagrams.
\item a circle that separates a jewel and a twisted band diagram. 
\end{enumerate}

We now describe the set of canonical Conway circles of a rational tangle.
\begin{definition}
  A {\bf rational tangle} of a projection $\Pi$ is a tangle $(\Delta ,\tau_{\Delta })$ such that all the {\bf canonical Conway circles} of $\Pi$ contained in $\Delta$ are concentric and delimit twisted annuli, with the exception of the innermost circle which contains $|a_m| \geq 2$ crossings of $\Pi$ in its inner disk (Figure 4).\\
 A {\bf maximal rational tangle} of a projection $\Pi$ is a rational tangle that is not strictly included in a larger rational tangle of $\Pi$.
 \end{definition} 
 \begin{figure}[h!]    
  \centering
   \includegraphics[scale=0.4]{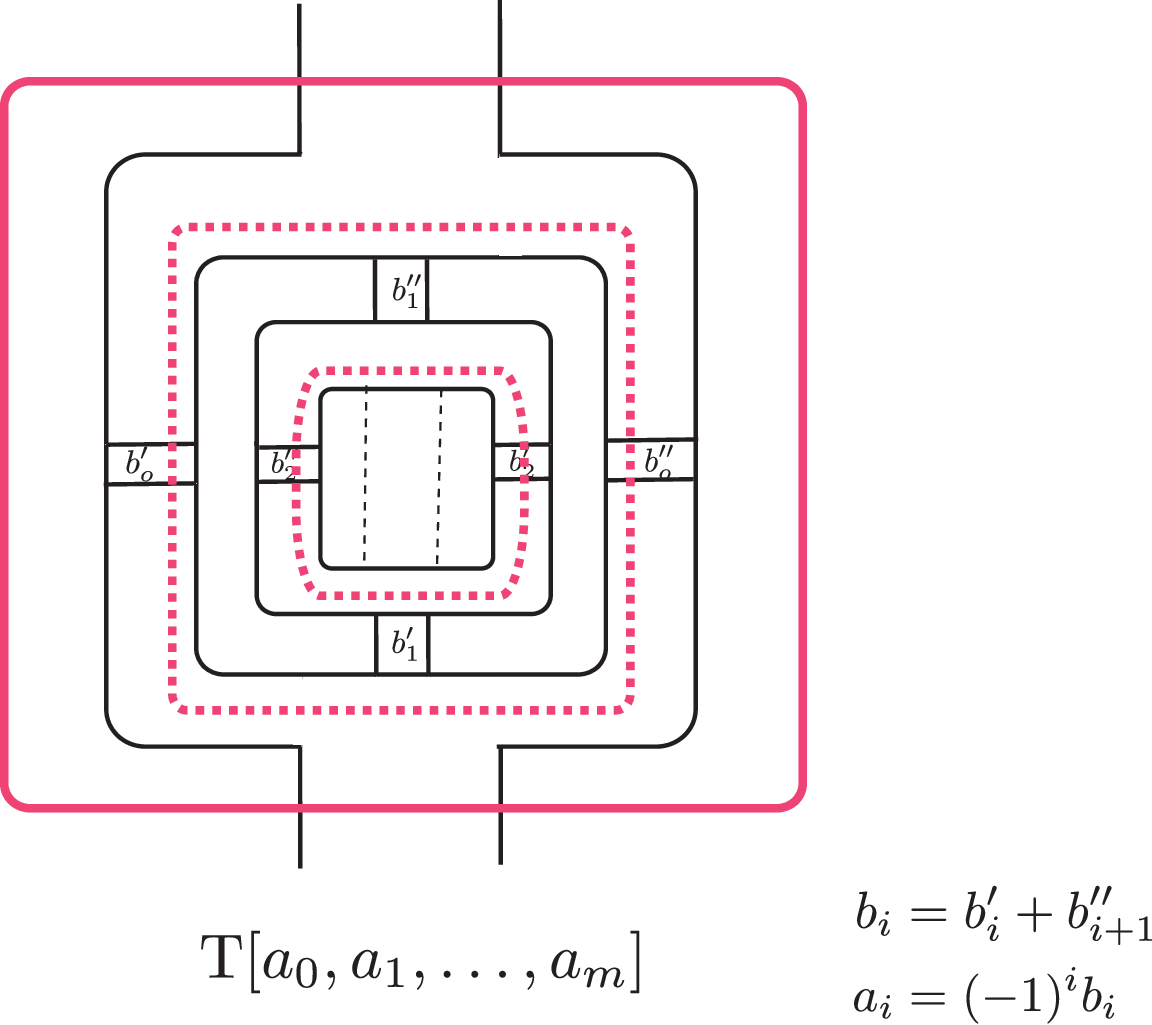}
\caption{${\rm T}[a_0, \dots, a_m]$}
\end{figure}

 Since we deal with alternating rational tangles, the weights $b_i$ have their sign alternate (two consecutive weights $b_k$ and $b_{k+1}$ have  $b_k b_{k+1} <0$). This implies that all the coefficients $a_i$ are of the same sign. \\
 Let $\rm T $ be a rational tangle  of $\Pi$ under its {\bf cardan} form  ${\rm T}[a_0, \dots, a_m]$ with $m \geq 1$ as shown in Fig. 4 and so that the first band  of weight $b_0$  is horizontal.\\
Since only alternating projections are involved, this implies that the signs of the weights $b_i$ alternate. This implies that the coefficients $a_i= (-1)^{i} b_i$ have all the same sign.\\
 We assign the continued fraction
   $$[a_0, a_1,  \cdots a_m]=a_0+ \dfrac{1}{a_1 +\dfrac{1}{{\ddots} +\dfrac{1}{a_m}}}$$
to the tangle $T$ under its cardan form ${\rm T}[a_0, \dots, a_m]$ where $a_0 \in  \mathbb{Z}$ and $a_1, \dots ,a_m \in \mathbb{Z}-\{0 \}$.
  
The rational number ${r \over s}=  [a_0, a_1, \cdots a_m]$ with $(r,s)=1$ and $r>0$ is called the {\bf fraction} of $\rm T$. The fraction is an isotopy invariant in the class of rational tangles.

Defining the {\bf numerator} of a tangle as the link obtained by joining the ``north" endpoints together and the ``south" endpoints also together, we can express any rational link $K$ as the numerator of an alternating rational tangle described in its cardan form of Fig.4. We denote ${\Pi}_{r \over s}$ this alternating projection and  $K_{r \over s}$ the corresponding knot.

\begin{remark}The set ${\cal C}_{can}({\Pi}_{r \over s})$ is non empty if and only if the fraction ${r \over s} \in {\mathbb Q} \setminus {\mathbb Z}$. The cardan form implies that the canonical Conway circles of ${\Pi}_{r \over s}$ are concentric. In the case of the projection ${\Pi}_{r \over s}$ with ${r \over s} \in {\mathbb Q} \setminus {\mathbb Z}$, we have the decomposition of $\Pi_{r \over s}$ into TBDs such that all are twisted annuli except two spires.
\end{remark}
\subsubsection{Essential Conway circles} \begin{definition} 
 An {\bf essential Conway circle} (abbreviated to {\bf essential circle}) of $\Pi$ is a canonical Conway circle that is not properly contained in a maximal rational tangle.
 \end{definition} 
 Let $\Pi$ be an (alternating) projection. Removing from the family ${\cal C}_{can}$ all concentric Conway circles of each maximal rational tangle of $(S^2, \Pi)$, we obtain the family of {\bf essential Conway circles} of $\Pi$ noted ${\cal C}_{ess}$. Therefore, the existence and uniqueness of ${\cal C}_{can}(\Pi)$ implies the existence and uniqueness of ${\cal C}_{ess}(\Pi)$. The decomposition of $(S^2,\Pi)$ by ${\cal C}_{ess}(\Pi)$ into TBDs, jewels and rational tangles is called the {\bf essential Conway decomposition} of $\Pi$. \\

For the rational projection $ {\Pi}={\Pi}_{r \over s}$ with ${r \over s} \in {\mathbb Q} \setminus {\mathbb Z}$, the set ${\cal C}_{ess}(\Pi)$ is empty.

\subsection{Flypes and Canonical Decomposition}

We can now locate precisely where flypes can be performed with respect to the canonical Conway decomposition of a reduced indecomposable alternating link projection.

\begin{theorem} \cite{quwe0}
(Position of flypes) Let $\Pi$ be a reduced indecomposable alternating link projection in $S^2$ and suppose that a flype can be performed in $\Pi$. Then, the active crossing point of the flype belongs to a diagram determined by $\mathcal{C}_{can}$ and such diagram is a TBD. The flype moves the active crossing point either within the twist region to which it belongs, or to another twist region of the same TBD.
\end{theorem}
\begin{definition}
 (1) The set of the twist regions of a given TBD is called a {\bf flype orbit}.\\
(2) The number of circles of $\mathcal{C}_{can}$ in the boundary of ${\mathcal T}$ is the {\bf valency} of the TBD $\mathcal T$.

\end{definition}
\begin{corollary} \cite{quwe0}\\
(1) A flype moves an active crossing point inside the flype orbit to which it belongs.\\
(2) Two distinct flype orbits are disjoint. 
\end{corollary}

This implies that an active crossing point belongs to a unique TBD. Since two TBDs have at most one canonical Conway circle in common, Corollary 3.1 can be interpreted as a loose kind of commutativity for the flypes.

The Flyping Theorem gives rise to an analogous theorem for tangles:

\begin{definition}
Let $T_1$ and $T_2$ be the alternating {\bf 2}-tangles with the same 4 endpoints as described by Fig.5. The transformation mapping $T_1$ on $T_2$ is called a flype on $T_1$ . \\
Two alternating {\bf 2}-tangles $\rm T$ and $\rm T'$ with the same 4 endpoints are {\bf flype-equivalent} if they are connected by a sequence of flypes on their {\bf subtangles} (tangles which are contained in $\rm T$ and $\rm T'$). They are then noted ${\rm T} \sim_f {\rm T'}$.

\end{definition}
\begin{remark}
The Flyping Theorem can be interpreted as follows:\\
- two alternating {\bf 2}-tangles with the same 4 endpoints are isotopic if and only if they are flype-equivalent and\\
- two reduced alternating projections $\Pi_1$ and $\Pi_2$ represent the same link isotopy class if they are connected by a sequence of flypes on their tangles, up to a homeomorphism on $S^2$.

 \end{remark}
  \begin{figure}[h!]    
  \centering
   \includegraphics[scale=.4]{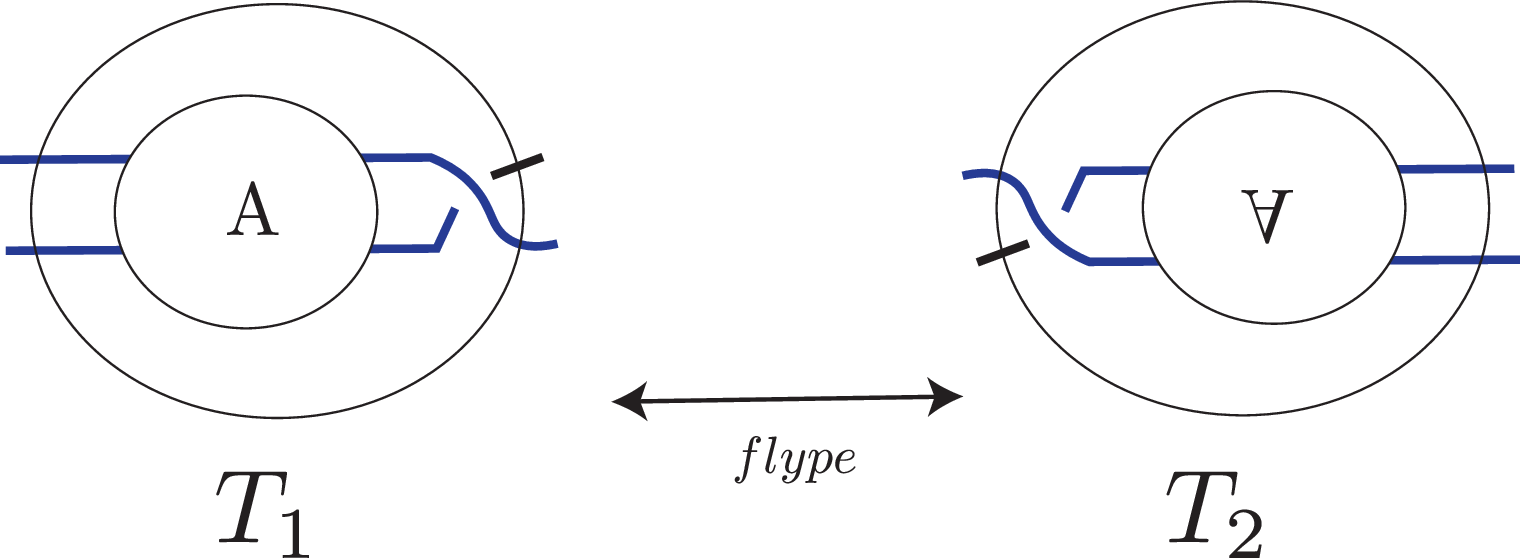}
\caption{$T_1$ and $T_2$ are related by a flype }
\end{figure}

 \begin{definition}
 Let $\rm T$ be a {\bf 2}-tangle whose extremities are NW, SW, SE and NE:\\
 Let $F$ be a rigid rotation of angle $\pi$ around an axis contained in the projection plane sending $\rm T$ into itself with NW on SW and SE on NE.
A {\bf 2}-tangle $\rm T$ is {\bf symmetric} if $\rm T \sim_f \rm F(T)$.
\end{definition}
\begin{example}  Fig.$6(a)$ shows a symmetric {\bf 2}-tangle: $\rm T \sim_f \rm F(T)$ and Fig.$6(b)$ shows a non-symmetric {\bf 2}-tangle $\rm T' \not \sim_f \rm F(T')$
\end{example}
 \begin{figure}[h!]    
  \centering
   \includegraphics[scale=.25]{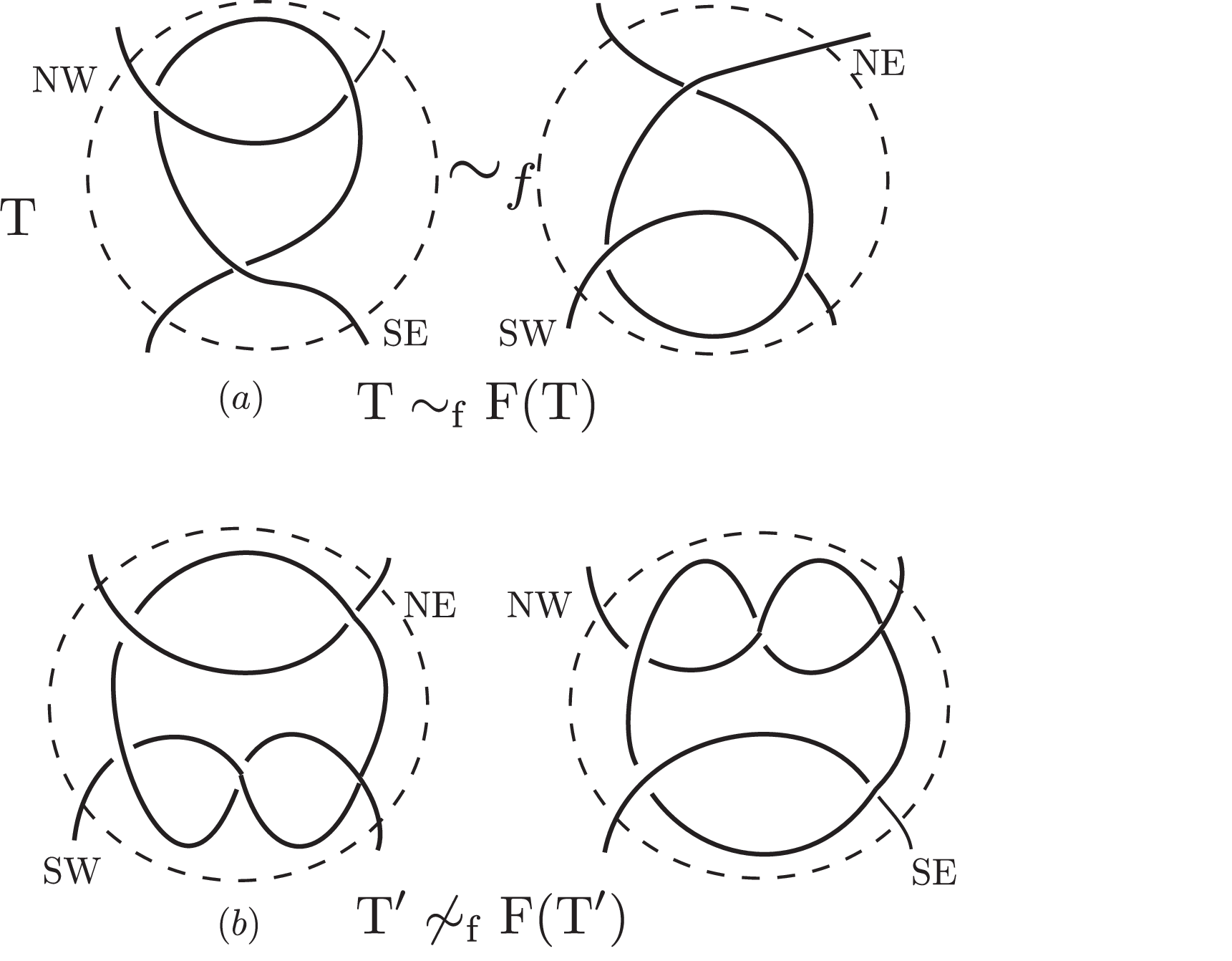}
\caption{$(a)$ A symmetric {\bf 2}-tangle and $(b)$ A non-symmetric {\bf 2}-tangle }
\end{figure}

\subsection{Canonical and Essential Structure Trees}

 Flyping Theorem is essential to deduce that two alternating projections in $S^2$ of the same isotopy class of a prime link in $S^3$ have their structure trees isomorphic. The structure trees are defined in \S 3.3.1 and \S 3.3.2.

Let $K$ be a prime alternating link and let $\Pi$ be an alternating projection of $K$.
\subsubsection{Canonical Structure Tree}
{\it Construction of the canonical structure tree ${\mathcal{A}} (K)$}.\\
Let $\mathcal{C}_{can}$ be the canonical Conway family for $\Pi$. The canonical structure tree ${\mathcal{A}} (K)$ is a graph
whose \\
-  vertices are in bijection with the diagrams determined by $\mathcal{C}_{can}$,\\
-  edges are in bijection with the circles of $\mathcal{C}_{can}$ so that the vertices of an edge associated to a circle $\gamma$ correspond to the diagrams which have $\gamma$ on their boundary.
The constructed graph ${\mathcal{A}} (K)$ is a tree.
 
We label the vertices of  ${\mathcal{A}} (K)$ as follows: if a vertex represents a TBD, we label it by its weight $a$ and if it represents a jewel, we label it with the letter $J$.

Since flypes and flat homeomorphism do not modify the canonical structure tree. Flyping Theorem implies that the canonical structure tree is independent of the choice of projection $\Pi$ provided it is alternating. Therefore, we can speak of the canonical structure tree of $K$ (and not of $\Pi$):
\begin{proposition} The canonical structure tree ${\mathcal{A}} (K)$ of a prime alternating link $K$ is independent of the alternating projection considered to represent $K$. 
\end{proposition}
\subsubsection{Essential Structure Tree}
{\it Construction of the essential structure tree $\mathcal\ \tilde A (K)$}.

On the same lines of the construction of the canonical structure tree ${\mathcal{A}} (K)$, we construct the {\bf essential structure tree} $\widetilde{\mathcal{A}}(K)$ as follows:\\
- its vertices are in bijection with the diagrams determined by the set $\mathcal{C}_{ess}(\Pi) $,\\
- its edges in bijection with the essential circles of $\mathcal{C}_{ess}(\Pi)$ and \\
 - the vertices of each edge $\gamma$ represent the two diagrams having the essential circle $\gamma$ in their boundary, \\ 
 - if a vertex represents a TBD, we label it by its weight $a$, if it represents a jewel, we label it with the letter $J$ and if it represents a rational tangle we label it with the corresponding fraction ${r \over s}$. 
  
As in the case of the canonical structure tree, we have:
  
\begin{proposition} The essential structure tree ${\mathcal{\widetilde A}} (K)$ is independent of the choice of a minimal projection of $K$. 
\end{proposition}
\vspace{1cm}

\begin{remark}
 The set $\mathcal{C}_{ess}(K)$  is empty if and only if $K$ has an (alternating) projection $\Pi$ such that $\Pi$ is either a jewel without boundary or a minimal projection of a rational knot/link. If $\mathcal{C}_{ess}(K)$  is empty, the essential structure tree $\widetilde{ \cal A} (K)$ is reduced to a single vertex. However, for a non-torus rational link $K_{r \over s}$, the canonical structure tree  ${\cal A} (K_{r \over s})$ is a linear tree such that the two end-vertices correspond to the two spires while the others correspond to the twisted annuli while the canonical structure tree of a jewel without boundary is a single vertex.
\end{remark}

\section{Lens Links and Free Periodic Links} 

By the Geometrization Theorem due to W. Thurston and G. Perelman, a $q$-freely periodic link in $S^3$ is topologically conjugate to a $(q,r)$-lens link for some integer $r$ such that $0 < r < q$ and $gcd (r,q)=1$. \\
In \S 4.2, we recall the description of the $(q,r)$-lens links as done in \cite{chbi}. In  \S 4.3, we make explicit the description of the $(q,r)$-visibility for torus knots of type $(n,m)$ and in particular for torus knots of type $(2,m)$ that are alternating.

\subsection {$N$-tangles}
We need a generalisation of Definition 3.4 of a two-dimensional tangle (with 2 strings) to a two-dimensional $N$-tangle (with $N \geq 2$ strings).
 \begin{definition}
 A two-dimensional $N$-tangle $\mathcal{T_N}$ is a pair $\mathcal{T}_N
=(\Delta ,\tau)$ where $\Delta $ is a disk and $\tau$ is a union of $N$ strings such that 
$\tau \cap \partial \Delta$ is composed of $2N$ endpoints with $N$ endpoints in the eastern half and $N$ endpoints in the western half.

A (two-dimensional) {\bf $N$-tangle} $\cal{T_N}$ of a projection $\Pi$  is a (two-dimensional) $N$-tangle such that $\Delta$ is a disk on the projection sphere $\mathbb{S}$ and $\tau= \Pi \cap \Delta$.
 \end{definition}
\begin{definition} A three-dimensional $N$-tangle is a pair $(B, \tau)$, where $B$
is a 3-ball and $\tau$ is a proper 1-submanifold of $B$ formed by $N$ strings whose $2N$ endpoints are in $\partial B$  such that $N$ are located in the eastern half of the equatorial circle and the others in its western half. 
Equivalently, a three-dimensional $N$-tangle is a pair $(\Delta \times [0,1], \tau)$ where $\Delta$ is a 2-disk and $\tau$ is a proper 1-submanifold of $\Delta \times [0,1]$ that meets $\partial(\Delta \times I)$ in $2N$ endpoints, half of which in $\Delta \times {0}$ and the other half in $\Delta \times {1}$.
 \end{definition}
 
A two-dimensional $N$-tangle is the projection of a three-dimensional $N$-tangle on the equatorial disk of $B$. 
\begin{remark}
 From a (two-dimensional) tangle $T$ of $\Pi$, we can create a three-dimensional tangle $\lambda (T)$  by a suitable small vertical perturbation near each crossing of $\Pi$. The ambient space of $\lambda (T)$ is taken to be a 3-ball for which the underlying disk of $T$ is an equatorial slice.
\end{remark}

 \begin{notation}
 Let $T= (\Delta ,\tau)$ be a tangle. If no misunderstanding is possible, we also call the 1-submanifold $\tau$ 
by $T$.
\end{notation}
We obtain a link projection from an $N$-tangle $\mathcal{T}$ by joining the $N$ right-handed strings to the $N$ left-handed strings by $N$ disjoint arcs on the projection sphere $\mathbb{S}$. The resultant link is called the {\bf closure of the $N$-tangle} $\mathcal{T}$.
 
 \begin{figure} 
\centering
\includegraphics[scale=.3]{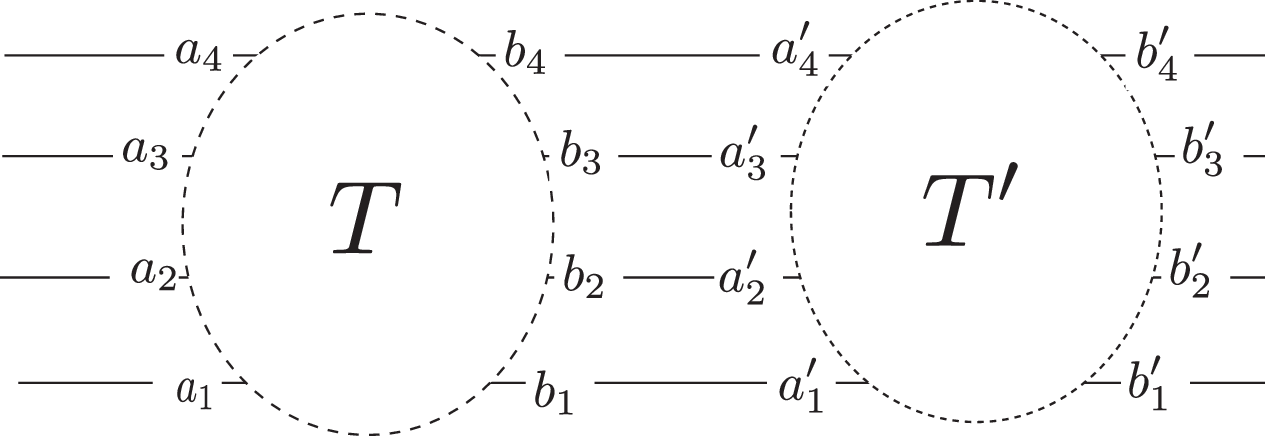}
\caption{ A tangle sum: $T \, T'$}
\end{figure}

Let $T$ and  $T'$ be two $N$-tangles contained in the projection sphere. We denote the $2N$ endpoints of  $T$ (respectively $T'$)  by $\{a_1,a_2 \dots, a_N \}$ (respectively $\{a'_1,a'_2 \dots, a'_N \}$) on its left side and $\{b_1,b_2 \dots, b_N \}$ (respectively $\{b'_1,b'_2 \dots, b'_N \}$) on its right side.\\
  We define the {\bf tangle sum} of $T$ and $T'$ by joining the points $\{b_1,b_2 \dots, b_N \}$ of $T$ to the points $\{a'_1,a'_2 \dots, a'_N \}$ of $T'$ with $N$ disjoint arcs contained in the projection sphere. The resulting tangle is denoted $T\,T'$ (Fig.7).

By using the terminology of the braid theory, let $\sigma_1, \sigma_2,\dots, \sigma_{N-1}$ be the $N$-tangles associated to the Artin generators of the braid group $B_N$. The element $\sigma_1 \sigma_2 \dots \sigma_{N-1}$ of $B_N$ gives rise to an $N$-tangle denoted 
\begin{equation} 
\Omega_N=\sigma_1 \sigma_2 \dots \sigma_{N-1}
\end{equation}
 Note that $ (\Omega_N)^N$ is the Garside element of the braid group $B_N$.
 
See Fig.8 for an example of a 4-tangle sum $ \rm T \,\rm T \, \rm T \, \Omega_4 \, \Omega_4 \, \Omega_4 ={\rm T}^3 (\Omega_4)^3$. 

\begin{figure}[h!]  
\centering
\includegraphics[scale=.4]{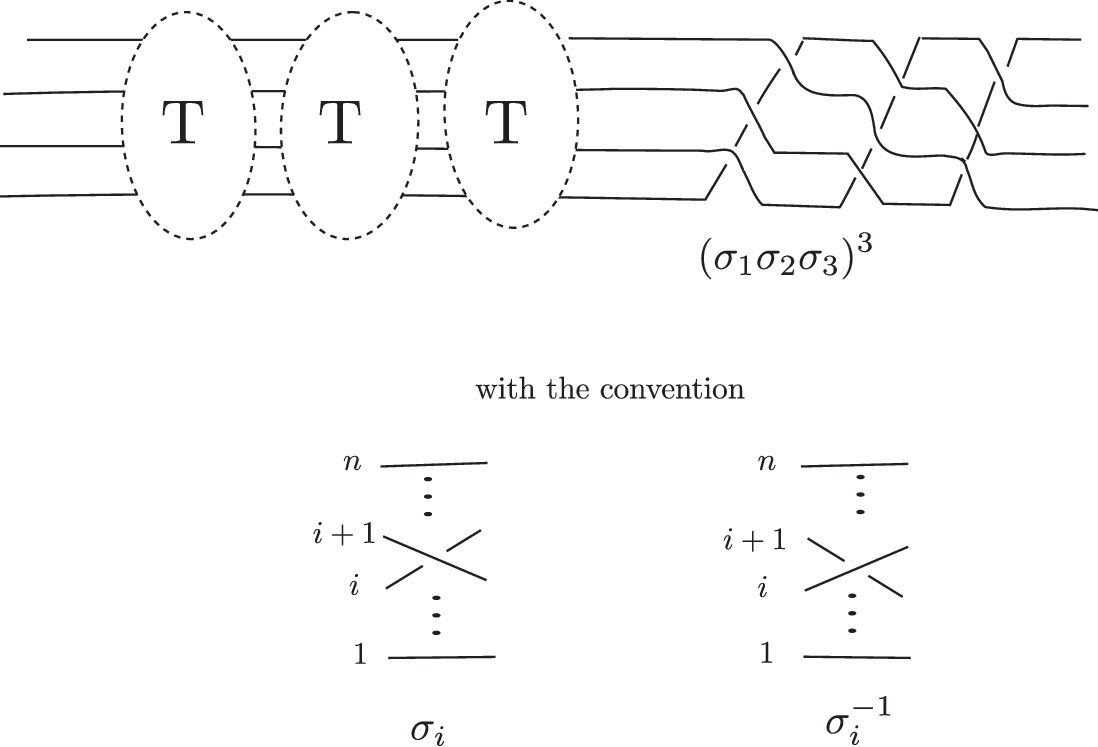}
\caption{The 4-tangle sum ${\rm T}^3 \, (\Omega_4)^3$ }
\end{figure}
  
\subsection{Description of the lens links and the $(q,r)$-visibility}
The quotient of $S^3$ by a free $q$-action of $\mathbb Z_q$ generated by a homeomorphism $\Phi$ onto $S^3$ of period $q$ without fixed point is a 3-manifold with fundamental group $\mathbb Z_q$. The Geometrization Theorem implies that the orbit manifold is a lens space $L(q,r)$ for some $r$ with $0< r < q$ and $gcd(r,q)=1$ and that $\Phi$ is topologically conjugate to the homeomorphism $\Phi_{q,r}$ defined by:
\begin{align}
\Phi_{q,r} &:S^3 \quad \longrightarrow  S^3  \nonumber \\
&(z_1,z_2)  \longmapsto  (\rho_q z_1, \rho_q^r z_2)
\end{align}
where $S^3=\{(z_1,z_2) \in \mathbb C \times \mathbb C : |z_1|^2 +|z_2|^2=1\}$ and
\begin{equation}
\rho_q=e^{i2\pi\over q} \quad {\rm and} \quad  \rho_q^r= e^{i2\pi r\over q}
\end{equation}

\begin{notation} The $q$-action generated by $\Phi_{q,r}$ is the {\bf $(q,r)$-action}. 
\end{notation}
 \begin{definition}
  Let  $q$ and $r$ be integers such that $0<r< q$ and $gcd(r,q)=1$.\\
  A link $L$ in $S^3$ is a {\bf $(q,r)$-lens link} if $L$ is invariant under the $(q,r)$-action. 
  \end{definition}
 \begin{remark} For $r=0$, the $(q,0)$-action is non-free. It is reduced to the rotation $R_q$ of angle $2 \pi \over q$ around the axis $ l_2=\{(z_1,z_2) \in S^3 |z_1=0\}$. A link $K$ disjoint from $l_2$ that is $\Phi_{q,0}$-invariant is $q$-periodic.
\end{remark}
The $(q,r)$-lens links are explicitly described as follows:
\begin{theorem}(\cite{chbi})
 Let $K$ be a (non-oriented) $(q,r)$-lens link. Then there exists an $N$-tangle $\rm T$ and the $N$-tangle $\Omega_N= \sigma_1\sigma_2 \dots \sigma_{N-1}$ (with $N \geq 2$) such that $K$ is the closure of  ${\rm T}^q (\Omega_N)^{Nr}$:
 \begin{equation}
 K= \widehat{{\rm T}^q (\sigma_1 \sigma_2 \dots \sigma_{N-1})^{Nr}}
 \end{equation}
 \end{theorem}
 See also \cite{man}
 
 \begin{definition}
 
The link $K$ expressed as the closure of ${\rm T}^q (\sigma_1 \sigma_2 \dots \sigma_{N-1})^{Nr}$ is said to be 
 {\bf $(q,r)$-visible} and $\rm T$ is called the {\bf fundamental lens-tangle}.
\end{definition}
 \begin{proof}  (Proof of Theorem 4.1)
 Let $K$ be a non-oriented $(q,r)$-lens link with $gcd(q,r)=1$ and $q \geq 2$. \\
 We express $S^3=\{ (z_1,z_2) \in {\mathbb C}^2 \, |\, |z_1|^2+|z_2|^2 =1 \}$, as the union of the two solid tori $V_1=\{(z_1,z_2) \in S^3 \, \, | \, |z_1| \geq |z_2| \} $ and $V_2=\{(z_1,z_2) \in S^3 \, \, | \, |z_1| \leq |z_2| \} $ with respectively the cores $l_1= \{(z_1,0 ) \in S^3 \}$ and $l_2= \{(0,z_2 ) \in S^3 \}$.\\
 We can assume that $K$ is in the interior $\mathring {V_1}$ of $V_1$. The action of $\mathbb Z_q$ in $\mathring {V_1}$  described by Equation (1) is the composition of $R_q$ and $S_r$:
 \begin{equation}
 \Phi_{q,r} (z_1,z_2)= S_r \circ R_q (z_1,z_2)= R_q \circ S_r (z_1,z_2) 
 \end{equation} 
where $R_q$ is the rotation around $l_2$ by an angle of $2\pi \over q$:
 \begin{equation}
 R_q(z_1,z_2)= (e^{i 2\pi \over q} z_1, z_2)
  \end{equation} 
and $S_r$ is the rotation around $l_1$ by an angle of $2\pi r \over q$:
 \begin{equation}
S_r(z_1,z_2)= (z_1, e^{i 2\pi r \over q} z_2).
 \end{equation} 
 \begin{figure}[h!]
\centering
\includegraphics[scale=.7]{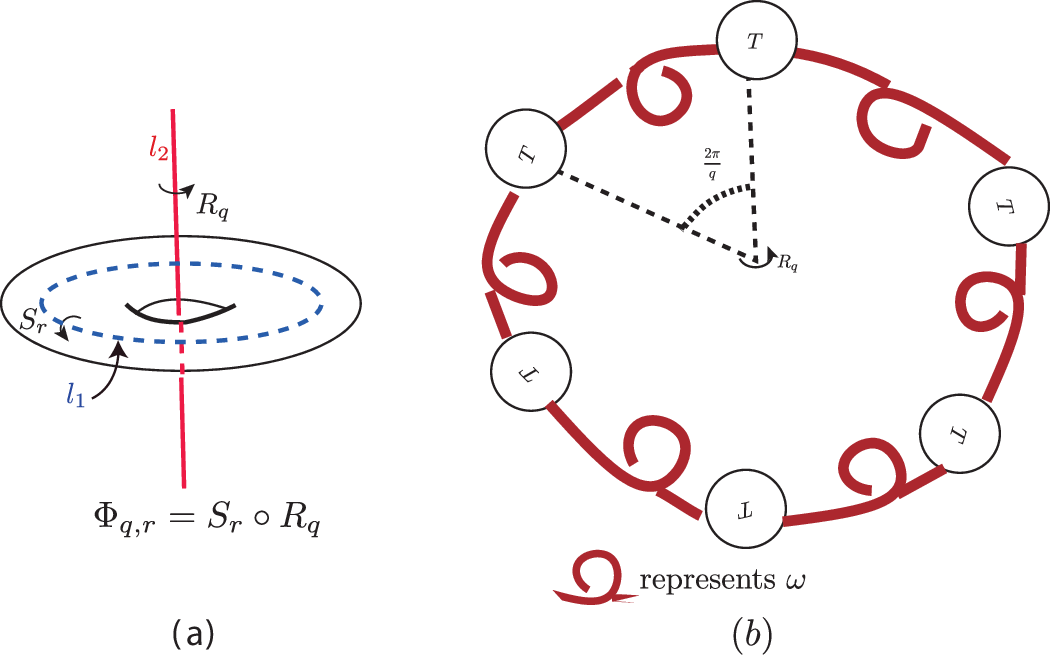}
\caption{$(a)$ Action of $\Phi _{q,r}$ \quad
 $(b)$ Schematic representation of a $(q,r) $-lens link }
\end{figure}

We now consider a fundamental domain $\cal D$ of the rotation $R_q$. There exists an $N$-tangle $\rm T$ for some $N$ and its corresponding three-dimensional $N$-tangle $\lambda (\rm T)$ (as constructed in Remark 4.1, $\lambda (\rm T)= K \cap \cal D$). Let $\rm O_N$ be the trivial $N$-braid and  $\lambda( \rm O _N)$ be its corresponding three-dimensional version.\\
The tangle $\lambda(\rm T)S_r(\lambda(\rm T))$ can be expressed as $\lambda(\rm T) \, \omega\lambda(\rm T)$ where $\omega$ denotes the 
 $N$-braid  $\lambda( \rm O _N)S_r(\lambda( \rm O _N))$.\\
If no misunderstanding is possible, to simplify, we also denote $\lambda(\rm T)$ by $\rm T$. Then $K$ is the closure of the tangle sum: 
\begin{equation} 
\rm T_1 \,\rm T_2\, \dots \,\rm T_q
 \end{equation} 
 where $\rm T_i= \rm T \, \omega$ for  each $i  \in \{1, \dots, q\}$. \\
Then $K$ can be expressed as the closure of 
 \begin{equation} 
 \rm T\, \omega \rm T \, \omega \dots \rm T \, \omega
  \end{equation} 
Since  $\rm T\, \omega =\omega\,  \rm T$, we obtain: 
\begin{equation}
K= \widehat {{\rm T}^q \, \omega ^q}.
\end{equation}
 We have
 \begin{equation}
 \omega^q=(\Omega_N)^{Nr}
\end{equation}
where $\Omega_N$ is given by Equation (1).\\
The $N$-tangle $(\Omega_N)^{Nr}$ denotes $(\rm O _N) S_r(\rm O _N) ... S_r(\rm O _N) = S_r^q (\rm O_N)$ where $S_r $ is a rotation of angle $2\pi r $ around $l_1$. 
So
\begin{equation}
\omega ^q= (\Omega_N)^{Nr}=(\sigma_1 \dots \sigma_{N-1})^{Nr}
\end{equation}
where $ (\Omega_N)^N$ denotes the braid resulting from the rotation of angle $2 \pi$ around $l_1$ of the $N$-braid $\rm O_N$.
\end{proof}

\begin{example}
\vspace{.3cm}
The 3-tangle ${\rm T}^q(\sigma_1 \sigma_2)^3$ in Fig.10 is $D_C( {\rm T}^q)$. 
\end{example}

\begin{figure}[h!]
\centering
\includegraphics[scale=.95]{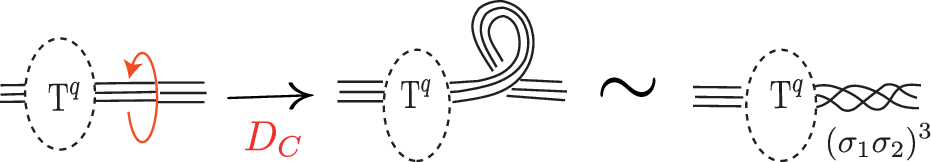}
\caption{ Right-handed Dehn twist along $C$ on ${\rm T}^q$}
\end{figure}

\begin{figure}[h!]
\centering
\includegraphics[scale=.5]{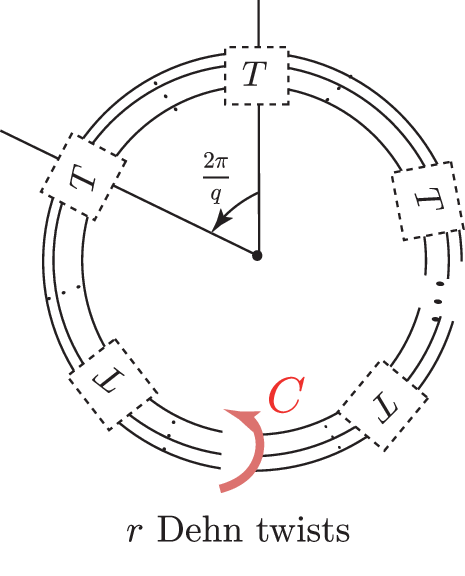}
\caption{$(q,r)$-visibility: $\widehat{{D_C}^r( {\rm T}^q)}$ }
\end{figure}

\begin{remark}
\begin{enumerate} 
\item Let $C$ be a circle bounding a sufficiently small disk transverse to the band connecting the $q$ tangles $\rm T$ as shown in Fig.10 and Fig.11. Performing the Dehn twist $D^r_C$ along $C$ (i.e. applying $r$ full twists to the band), we obtain
$D_C^r(\widehat{{ \rm T}^q})= \widehat{D^r_C( {\rm T}^q)}$.\\
For $r>0$, these twists are right-handed. \\
The operation $D_C^r$ transforms the link $\widehat{{ \rm T}^q}$ into the $(q,r)$-lens link $\widehat {{\rm T}^q(\sigma_1 \dots \sigma_{N-1})^{Nr}}$. \\
By identifying $l_1$ with $C$, we can identify $D^r(C)$ with $S_r^q$ where $S_r$ is described by Equation (7).\\

\item Given a $(q,r)$-lens link described by Equation (10), we visualize its $q$-symmetry in the 3-space by distributing 
 $q$ times of the $1 \over q$-th of $r$ twists along the band that corresponds to $\omega$.  
\end{enumerate}
\end{remark}

 \subsection {On the free periodicity of the alternating torus knots and their $(q,r)$-visibility.} 
 We first recall a result of R.Hartley on a necessary and sufficient condition of a torus knot of type $(m, n)$ to be $q$-freely periodic.
 \begin{lemma}(\cite{hartl})
  A torus knot of type $(m,n)$ has free period $q$ if and only if $q$ is coprime with $m n$.
  \end{lemma}

  To be a $(q,r)$-lens link where $0< r<q$, the following condition is necessary and sufficient:
\begin{lemma}
 Let ${\rm T}(n,m)$ be a torus link of type $(n,m)$. Then ${\rm T}(n,m)$ is a $(q,r)$-lens link if and only if $r n-m$ or $rm-n$ is divisible by $q$.
\end{lemma}
\begin{proof}
The link ${\rm T}(n,m)$ can be seen as the intersection of $S^3=\{ (z_1,z_2) \in {\mathbb C}^2 \, \, | \, |z_1|^2+|z_2|^2 =1 \}$ with the complex surface $\Sigma$ defined by :
\begin{equation}
 \Sigma = \{(z_1,z_2) \in {\mathbb C}^2 \, | \,z_1^m+z_2^n=0\}
\end{equation}
Since $$\Phi_{q,r}(z_1, z_2)= (\rho_q z_1, \rho_q^{r} z_2),$$ we have
$$\rho_q^m z_1^m+ \rho_q^{nr} z_2^n=\rho_q^m(z_1^m+\rho_q^{nr-m} z_2^n).
 $$  
Therefore, the torus link ${\rm T}(n,m)$ link is a $(q,r)$-link if and only if 
\begin{equation}
r\,n-m \, \, \text{is divisible by} \, \,q.
\end{equation}
\end{proof}

Assume that $m \neq \pm 1$. A torus knot is alternating if and only if it is a knot ${\rm T}(2,m)$
 with $m \equiv 1\mod{2}$.

{\bf Question:} For what pair of integers $(q,r)$, is a torus knot of type $(2,m)$ a $(q,r)$-lens knot?

By Lemma 4.1(\cite{hartl}), if $q$ is coprime with $2m$, $q$ is a free period of ${\rm T}(2,m)$.

Therefore, for each alternating knot ${\rm T}(2, m)$, there is an infinite number of free periods $q$ for every $q$ coprime to $2m$.

Assume that $m \geq 3$. To obtain the form described by Equation (10), it is sufficient to solve the following diophantine equation with unknowns $r$ and $k$.
\begin{equation}
2{\bf r}- q{\bf k} = m
\end{equation}

 Therefore, if $(r,k)$ is a solution of Equation (15) with $0<r<q$, ${\rm T}(2,m)$ is a $(q,r)$-lens knot. 
\begin{example}
Consider the trefoil knot ${\rm T}(2,m=3)$.  By Lemma 4.1, since $gcd(2,m)=1$, every $q$ coprime with $2m=6$ is a free period of the knot. Hence $q$ is odd and $q>3$. \\
Since $m=3$, Equation (15) has the form:
\begin{equation}
 2{\bf r}-q{\bf k} = 3
 \end{equation}

For instance, ${\rm T}(2,3)$ is freely 5-periodic. The pair $(r,k)=(4,1)$ is a solution of Equation (16) and ${\rm T}(2,3)$ is a  $(5,4)$-lens knot.

\end{example}

\section{Visibility of $q$-actions} 

This section deals with a generalization of the Visibility Theorem 1.1 to the case of free $q$-actions. 

\subsection{Flypes and Flyping Theorem} 

Let $\Pi$ be an $n$-crossing regular projection of a knot $K$ on the projection sphere.  As in \cite{meth}, consider $n$ disjoint small ``crossing balls " $B_1,\dots , B_n $ neighbourhoods  of the crossing points $c_1, \dots c_n$ of $\Pi$. Then, assume that $K$ coincides with $\Pi$, except that inside each $B_i$ the two arcs forming $\alpha(c_i) =\Pi \cap B_i$ are perturbed vertically to form semicircular overcrossing and undercrossing arcs which lie on the boundary of $B_i$. This representation of the knot $K$ from a projection $\Pi$ is expressed as $K= \lambda(\Pi)$ (\cite{meth}). 
We call $\lambda (\Pi)$ a {\bf realized projection} (or a realized diagram in the terminology of \cite{boy} of $K$).

Considering the ambient space $S^3$ as ${\mathbb R^3} \cup\{\infty\}$, we shall take the projection 2-sphere $\mathbb{S}$  to be $ S^2=\{\bold x \in \mathbb R^3 :\bold ||x||=1\}$. We assume that $\lambda(\Pi)$ lies within the neighborhood $N= \{x \in \mathbb R^3: {1\over 2} \leq ||{\bold x}  || \leq {3 \over 2} \}$.

\begin{definition} Let $g : (S^3, \lambda (\Pi_1)) \rightarrow  (S^3, \lambda (\Pi_2))$ be a homeomorphism of pairs. The homeomorphism of pairs $g$ is {\bf flat} if $g$ is isotopic to a homeomorphism of pairs $h$ with the condition that $h$ maps $N$ onto itself
and $h|_N = h_0 \times  id_{[{1\over 2}, {3\over2}]}$ for some orientation-preserving homeomorphism $h_0 : S^2 \rightarrow S^2$. We call $h_0$ the {\bf principal part of the flat homeomorphism} $g$.
\end{definition}

\begin{definition} An {\bf isomorphism of realized projections} $h: (S^3, \lambda (\Pi)) \rightarrow  (S^3,\lambda ({\widetilde \Pi}))$ is a homeomorphism of pairs $h:(S^3, K) \rightarrow (S^3, \widetilde K)$ such that:\\
(1) $h(S^2)=S^2$,\\
(2)  $h(B_i) ={\widetilde B_i}$\\
(3)  $h(\alpha(c_i))=\alpha(\widetilde c_i)$.

\end{definition}
A flat homeomorphism of pairs $g$ as defined in Definition 5.1 is an isomorphism of realized projections.

We recall the definition of a flype as described in \cite{meth}.

\begin{definition} Let $\Pi_1$ be a projection with the pattern described in Fig.$3(a)$.
 A {\bf standard flype} of $f_s:(S^3, \lambda(\Pi_1)) \rightarrow (S^3, \lambda(\Pi_2))$ where $\Pi_2$ is the pattern depicted in Fig.$3(b)$, is a homeomorphism of pairs such that:\\
 (1) $ f_s$ sends the 3-ball $B_{\rm A}$ into itself by a rigid rotation around an axis in the projection plane,\\
(2) $ f_s$ fixes the 3-ball $B_{\rm B}$ pointwise,\\
(3) $ f_s$ moves the visible crossing on the left of Fig.$3(a)$ to the visible crossing on the right of Fig.$3(b)$.
\end{definition}

\begin{definition}  
Let $\Pi_1$ be any projection. Then a flype is any homeomorphism of pairs
$f: (S^3, \lambda(\Pi_1))  \rightarrow  (S^3, \lambda(\Pi_2))$ of the form $ f= g_1 \circ  f_s \circ  g_2$ where $ f_s$ is a standard flype and $g_1$ and $g_2$ are flat homeomorphisms.
 \end{definition}
 
In the case where the tangle $A$ contains no crossing, the standard flype described in Fig.3  is a flat homeomorphism. Therefore, according to the above definition, any flat homeomorphism can be considered as a flype and we call it a {\bf trivial flype}. 
We reformulate the Menasco-Thislethwaite Flyping Theorem as follows:

\begin{theorem} (Flyping Theorem). Let $\Phi:(S^3,K) \rightarrow (S^3,K)$ be any orientation-preserving homeomorphism of pairs where $K$ is a prime alternating knot. Let $\lambda(\Pi)$ be a realized projection of a reduced alternating projection $\Pi$ of $K$ and $h:(S^3,K)  \rightarrow (S^3,\lambda (\Pi))$ be a homeomorphism of pairs. Then the isomorphism of realized projections $\Phi_\Pi$ $(= h \circ \Phi \circ h^{-1})$ of $(S^3,\lambda (\Pi))$ can be expressed as
$$
\Phi_\Pi =\phi \circ F
$$
where $\phi$ is a flat homeomorphism and $F$ is a composition of standard flypes on $\lambda (\Pi)$ unless $F$ is the identity map.

\end{theorem}
\begin{remark}
\begin{enumerate}
\item
The standard flypes involved are not trivial.
\item
Standard flypes and flat homeomorphisms  ``essentially " commute  (by ``essentially "  commute, we mean that if $f_s$ is a standard flype and $h$ is a flat homeomorphism, then $f_s \circ h = h \circ f'$ where $f'= (h^{-1} \circ f _s \circ h)$ is also a standard flype. This is the reason why it is possible to express $\Phi _{\Pi }=\phi \circ F$.
\end{enumerate}
\end{remark}

For a $q$-action, we have Flyping Theorem 5.1 under the following form.
\begin{theorem} (Flyping Theorem for a $q$-action).
 Let $K$ be a prime oriented alternating knot with a $q$-action ($q \geq 3$) which is defined by a $q$-homeomorphism of pairs $\Phi: (S^3, K) \rightarrow (S^3, K)$. Let $\lambda(\Pi)$ be a realized projection of a reduced alternating projection $\Pi$ of $K$ and $h:(S^3,K)  \rightarrow (S^3, \lambda (\Pi))$ be a homeomorphism of pairs. Then the isomorphism of the realized projection $(S^3, \lambda (\Pi))$ onto itself $\Phi_\Pi (= h \circ \Phi \circ h^{-1})$ can be expressed as:
 \begin{equation} 
\Phi_\Pi =\phi \circ F
\end{equation}
where $\phi$ is a flat homeomorphism of order $q$ and $F$ is a composition of standard flypes on $\lambda (\Pi)$ unless $F$ is the identity map.\\
The principal part $\phi ^*$ of $\phi $ on the projection sphere $\mathbb S$ is conjugate to a rotation of angle $2 \pi \over q$ without fixed point on $\Pi$.
\end{theorem}
\begin{proof} 
Since $\Phi_\Pi= \phi \circ F$ satisfies $\Phi_\Pi^q ={\rm  Id}$, by Remark 5.1.2 we have $\Phi_\Pi^q=\phi^q \circ {\tilde F} = {\rm  Id}$. This implies that the restriction $\phi ^*$ of $\phi$ on $\mathbb S$, can be constructed as an automorphism of order $q$ of $\mathbb S$. By Kerekjarto's theorem (\cite{coko}), $\phi ^*$ on $\mathbb S$ is conjugate to a rotation of angle $2 \pi \over q$ without fixed point on $\Pi$.
 \end{proof}
 
 \begin{notation}
 By abuse of notation, we also denote the principal part $\phi^*$ by $\phi$. 
\end{notation}

By Flyping Theorem, we have two cases to treat: 
  
(1) Case where there are no flypes. Thus, $\Phi_\Pi$ is conjugate to a $q$-rotation and the $q$-action is necessarily non-free. 

(2) Case where flypes are involved. This happens for some non-free $q$-actions but it is always the case for any free $q$-action, as we will see later.

\subsection{Statement of Visibility Theorems for a $q$-action} 

Let $K$ be a prime alternating knot in $S^3$ that is not a rational knot. Let $\Pi$ be a reduced alternating projection and $\lambda(\Pi)$ a realized projection of $\Pi$ and  $h:(S^3,K)  \rightarrow (S^3,\lambda (\Pi))$ be a homeomorphism of pairs. Let $\Phi:(S^3,K) \rightarrow (S^3,K)$ be an orientation preserving homeomorphism of pairs, thus $\Phi_\Pi = h \circ \Phi \circ h^{-1}$ is an isomorphism of the realized projection $\lambda (\Pi)$ onto itself.
 
In the case where $\Phi$ is $q$-periodic, the Visibility Theorem 1.1 is equivalent to the following theorem:
 \begin{theorem} (Visibility Theorem for a non-free $q$-action).
 Let $K$ be a prime oriented alternating knot with a non-free $q$-action ($q \geq 3$) defined by a $q$-periodic homeomorphism of pairs $\Phi:(S^3,K) \rightarrow (S^3,K)$.
 Then there exists a reduced alternating projection $\Pi$ of $K$ such that $\Phi$ is isotopic via pair maps to an isomorphism $\Phi_\Pi$ of a realized projection $(S^3, \lambda (\Pi))$ onto itself which is reduced to a flat homeomorphism $\phi$ of order $q$:
\begin{equation}
 \Phi_\Pi =\phi 
\end{equation}
The restriction $\phi ^*$ of $\phi$ on the projection sphere $\mathbb S$ is a rotation of angle $2 \pi \over q$ without fixed points on $\Pi$.
\end{theorem}

In the case where $\Phi$ is freely $q$-periodic, the composition of flypes $F$  in $\Phi_\Pi =\phi \circ F$ (Equation (17))  is essential to describe the isomorphism $ \Phi_\Pi$. The composition of flypes $F$ cannot be reduced to the identity map because if this were the case, $\Phi_\Pi$ would be conjugate to a rotation of order $q$ and this would contradict the assumption that the $q$-map $\Phi$ is fixed point free. Specifically, we have:

 \begin{theorem} (Visibility Theorem for a free $q$-action). 
 Let $q$ be $ \geq 3$. Let $K$ be a prime oriented alternating knot that is with a $q$-action defined by a   
 free $q$-periodic homeomorphism of pairs $\Phi:(S^3,K) \rightarrow (S^3,K)$.
 Then there exists an alternating projection $\Pi$ of $K$ with a twisted band diagram $\mathcal T_0$ such that $\Phi$ is isotopic to an isomorphism $ \Phi_\Pi$ of of a realized projection $(S^3, \lambda (\Pi))$ onto itself that satisfies the following conditions:
\begin{equation}
 \Phi_\Pi =\phi \circ F_0
 \end{equation}
where \\
(1) $\phi$ is a flat homeomorphism of order $q$ whose principal part is topologically conjugate to a rotation, leaving invariant the diagram $\mathcal T_0$ and with generic orbits for its action on the other diagrams of the essential Conway decomposition, and\\
(2) $F_0$ is a non-trivial composition of standard flypes all performing on $\lambda(\mathcal T_0)$.
 \end{theorem}
The TBD $\mathcal T_0$ is called the {\bf main TBD} of $\Pi$. 
\begin{remark}
The above visibility theorems provide necessary conditions on the essential structure tree for $K$ to admit a free or non-free $q$-action. For example, for an alternating knot to be freely $q$-periodic, the essential tree must admit an automorphism of order $q$ such that its fixed set is a vertex corresponding to a TBD (which is the main TBD).
\end{remark}
 
 \subsection{Action of $\phi$ on the essential decomposition of $(S^2, \Pi)$}  

In the following, an {\bf automorphism} of $S^2$ preserves the orientation of $S^2$.

Suppose that we have a finite decomposition of the 2-sphere $S^2$ ($\subset S^3)$ into connected planar surfaces $S_k$ such that
$S_i \cap S_j$ is empty or is a boundary component for $i \neq j$.

In addition, let us assume that we have an automorphism $g : S^2 \rightarrow S^2$ of order $q$ which respects the decomposition, that is to say that for every index $i$, there exists an index $k(i)$ such that $g(S_i) = S_{k(i)}$. By Kerekjarto's theorem \cite{coko}, $g$ is conjugate to a rotation of order $q$. Let $\gamma_k$ be a boundary component of the surfaces $S_i$ and $S_j$ of the decomposition. Then we have the two following possible cases:
\begin{enumerate}
\item $g(S_i) , g^2(S_i) , \dots , g^q(S_i) = S_i$ are all distinct. We say that the orbit of $S_i$ under $g$ is {\it generic} (or more simply $S_i$ is {\bf generic}).  

\item $g(S_i)=S_i$. Then $g|_{S_i}$ is an automorphism of $S_i$. We say that $S_i$ is $g$-invariant. .
\end{enumerate} 
 Let $\gamma_k$ be a boundary component of a surface $S_i$ of the decomposition. Then we have the two following cases:
\begin{enumerate}
\item $g(\gamma_k) , g^2(\gamma_k) , \dots , g^q(\gamma_k) = \gamma_k$ are all distinct. We say that the orbit $\gamma_k$ under $g$ is {\it generic} (or $\gamma_k$ is {\bf generic}).  
\item $\gamma_k$ is $g$-invariant. 
\end{enumerate} 
(see {\cite{erquwe}).
As above, let $K$ be a (prime oriented) alternating knot in $S^3$ with a $q$-action defined by a $q$-homeomorphism of pairs $\Phi: (S^3,K) \rightarrow (S^3,K)$ and let $\Pi$ be a reduced alternating projection of $K$. The restriction of the homeomorphism $\phi$ on $(S^2, \Pi)$ expressed as in Equation (17) is a periodic automorphism of $S^2$. Each component $S_i$ given by the essential (or canonical) Conway decomposition is generic or invariant by $\phi$.

Let ${\cal T}_i$ be a {\bf 2}-tangle ${\cal T}_i=(\Delta_i, \tau_i)$ such that the boundary  $\partial {\cal T}_i$ is an essential (or canonical) Conway circle $\gamma_i$.
By Kerekjarto Theorem, $\phi$ is topologically conjugate to a rotation then each essential (or canonical) circle $\gamma_i$ is invariant or generic.
 
\begin{proposition} 
Suppose that an essential  (or canonical) Conway circle $\gamma_i$ is $\phi$-invariant. Then $q=2$.

\end{proposition}
\begin{proof}
 By the property that a {\bf 2}-tangle has two entry points and two exit points, the only possible non-trivial automorphism implies that $q=2$.
\end{proof} 
 \begin{proposition}
 If $q \geq 3$, any essential  (or canonical) Conway circle is generic.
 \end{proposition}
 \begin{proof}
 It is a straightforward consequence of Proposition 5.1.
  \end{proof}

  Since the canonical structure tree ${\cal A} (K_{r \over s})$ is linear, we have the following restrictions on $q$ for a $q$-action on a non-torus rational knot $K_{r \over s}$.
  \begin{proposition} 
  Let $K_{r \over s}$ be a non-torus rational knot. Then there is no $q$-action with $q \geq 3$.
  \end{proposition}
\begin{proof}
We will focus on the set ${\cal C}_{can} (K_{r \over s})$ and the canonical structure tree ${\cal A} (K_{r \over s})$ which is linear by Remark 3.3.\\
Suppose the opposite: let $K_{r \over s}$ be a non-torus rational knot endowed with a $q \geq 3$-action $\Phi$. By Proposition 5.1, there is no $\phi$-invariant canonical circle. This implies that every circle of ${\cal C}_{can} (K_{r \over s})$ is generic. Since the edges of ${\cal A} (K_{r \over s})$ correspond to the circles of ${\cal C}_{can} (K_{r \over s})$, the orbit of each edge under $\varPhi$ is generic. Moreover, since ${\cal A} (K_{r \over s})$ is a linear tree, each non-trivial automorphism on ${\cal A} (K_{r \over s})$ is of order 2. This contradicts the fact that the order of the automorphism $\varPhi$ on ${\cal A} (K_{r \over s})$ induced by $\Phi$ should be $q \geq 3$.

\end{proof}

 \subsection{Action of $\tilde{\varPhi}$ on Structure Trees} 
In this section, we will study the action of $\tilde{\varPhi}$ on the essential structure tree $\widetilde{\mathcal{A}}(K)$. We will see that the set $ Fix (\widetilde \varPhi)$ is essential to give informations on the visibility of the $q$-action. 
 
\subsubsection{Analysis of $ Fix (\widetilde {\varPhi})$}  

For the proof of the Visibility Theorems of the $q$-actions, we need a detailed analysis of the action of the automorphism $\widetilde{\varPhi}$ induced by $\Phi_\Pi$ on the essential structure tree $\widetilde{ \cal A} (K)$.
\begin{proposition}  
The automorphism $\widetilde{\varPhi}$ on $\widetilde{\mathcal{A}}(K)$ satisfies $\widetilde{\varPhi}^q={\rm Identity}$. 
\end{proposition} 
\begin{proof} It is a direct consequence of $\Phi_\Pi^q= \rm Identity$.
\end{proof}

Assume that $K$ is an oriented alternating knot endowed with a $q$-action. 

We first consider the case where ${\widetilde A}(K)$ is reduced to a single vertex. Then $K$ is either a rational knot or a jewel without boundary.\\
$(a_1)$ {\bf $K$ is a rational torus knot}. Then $K$ is a torus knot of type $(2, m)$ with $m$ an odd integer $\geq 3$. This implies that a torus knot ${\rm T}(2,m)$ has periods 2 and $m$. By Lemma (4.1), ${\rm T}(2,q)$ has infinitely many free periods for any $q$ coprime to $2m$.\\
$(a_2)$ {\bf $K$ is a non-torus rational knot}. A non-torus rational knot $K_{r \over s}$ cannot have a $q$-action with $q \geq 3$ (Proposition 5.3).

$(b)$ {\bf $K$ is a jewel without boundary}. Since there are no possible flypes on any alternating projection of $K$, the $q$-action is necessarily equivalent to a flat homeomorphism of period $q$ with principal part topologically conjugate to a rotation. Thus, the $q$-action is non-free.

From now, we assume that ${\widetilde A}(K)$ is {\bf not reduced to a single vertex}.

Since the graph $\widetilde{\mathcal{A}}(K)$ is a tree, the fixed point set $ Fix (\widetilde \varPhi)$ is a non-empty subtree. We then have two cases to deal with:
\begin{enumerate}
\item $Fix (\widetilde \varPhi)$ contains an edge $E$. 
\item $Fix (\widetilde \varPhi)$ is reduced to a single vertex $V_0$.
\end{enumerate}

We now describe these two cases of $Fix (\widetilde \varPhi)$ in terms of the essential decomposition of $(S^2,\Pi)$. 

\subsubsection{Case where there is an edge $E$ invariant by $\widetilde \varPhi$}

The case of an edge invariant by $\widetilde {\varPhi}$  corresponds to the case of the existence of an essential circle invariant under $\phi$. This implies that $q=2$ by Proposition 5.2.\\
\begin{proposition} If $\Phi$ is of order $q >2$, $Fix (\widetilde \varPhi)$ is reduced to a a vertex $V_0$ of $\widetilde {\cal A}(K)$.
\end{proposition} 
\begin{proof}
 If $q \geq 3$, by Proposition 5.2, the $\phi$-orbit of each essential circle is generic which means that there is no $\phi$-invariant essential circle. Therefore, there would be no invariant edge under $\widetilde \varPhi$ and $Fix (\widetilde \varPhi)$ is reduced to a single vertex of $\widetilde {\cal A}(K)$.
\end{proof}

\subsection{Proof of Visibility Theorems} 
For the visibility of a $q$-action, let us highlight the existence or absence of an invariant projection 2-sphere under the $q$-action involved. The Flyping Theorem 5.2 gives rise to the situation described by Equation (14):
$$ \Phi_\Pi= \phi \circ F$$
where $\phi$ is of order $q$ and $F$ is a composition of flypes.

$ \blackdiamond$ For a non-free action $\Phi$ (\S 2), there exists a circle $\alpha \in S^3$ which is fixed by $\Phi$ and which cuts the projection 2-sphere $\mathbb S$ in two points. Thus, $\mathbb S$  is invariant under $\Phi_\Pi$ and $\phi$. (Note that each $q$-periodic knot can be identified with an invariant knot under $\Phi_{q,0}=S_r\circ R_q$ of Equation (4)) and the axis $\alpha$ can be identified with the rotation axis $l_2$ of $R_q$.)

$ \blackdiamond \blackdiamond$ For a free $q$-action, the isomorphism $\Phi_\Pi$ has no fixed point on the sphere $S^3$.  For freely $q$-periodic links described as invariant links under $\Phi_{q,r}$ with $0 <r < q$ and $gcd(r,q)=1$, the homeomorphism $S_r$ cannot be trivial and $\mathbb S$ is not invariant  under $\Phi_{q,r}$. This implies that the composition of flypes $F$ in Equation (13) cannot be reduced to the identity map no matter what alternating projection $\Pi$ is considered.  

Since $q >2$, Proposition 5.5 implies that $Fix (\widetilde {\varPhi})$ is reduced to a single vertex $V_0$ which corresponds either to a jewel $J_0$ or to a TBD $\mathcal T_0$. We have two cases to analyze:
\begin{enumerate}

\item {\bf Case 1}, $V_0$ corresponds to a jewel $J_0=(\Sigma_0, \Sigma_0 \cap \Pi)$ such that its boundary  is a set of $k$ essential circles $\gamma_1, \dots, \gamma_k$. Each essential circle $\gamma_i$ is  the boundary of a tangle ${\cal T}_i=(\Delta_i,\tau_i= \Pi \cap \Delta_i)$ whose disk $\Delta_i$ does not meet the interior $\mathring J_0$.\\ 
Since $J_0$ is a jewel, no flypes can occur on it. When $Fix (\widetilde \varPhi)$ corresponds to a jewel, we have Visibility Theorem 3.1, as done in \cite{co2}. 

 \begin{proposition} 
If $Fix (\widetilde \varPhi)$ corresponds to a jewel, the knot $K$ is necessarily $q$-periodic and it admits a $q$-periodic alternating projection.
\end{proposition}
Then if $Fix (\widetilde \varPhi)$ corresponds to a jewel, we obtain Theorem 5.3.
We then have an interesting relation between the essential decomposition and the free periodicity for the class of alternating knots:
\begin{proposition} 
If $K$ is freely $q$-periodic, $Fix (\widetilde \varPhi)$ corresponds to a twisted band diagram.
\end{proposition} 
\item {\bf Case 2}. $V_0$ corresponds to a TBD $\mathcal{T}_0$. Let $V_0$ be of valency $\mu$ and of weight $w$.\\
The case where $\mu=0$ and $w$ is an odd integer, corresponds to a torus knot of type $(2,w)$.\\ 
 Assume that $|w| \neq 1$. Then $V_0$ corresponds to a TBD ${\cal T}_0$ - the {\bf main TBD} of $\Pi$- which is $\Phi_\Pi$-invariant.
 The valency $\mu$ of $V_0$ is a multiple of $q$ since the $\phi$-orbit of each essential circle is generic. 
Let $\gamma^*_1, \dots \gamma^*_{\mu}$ be the $\mu$ essential circles which constitute the boundary of ${\cal T}_0$. The isomorphism $\Phi_\Pi$ induces an automorphism on these circles.
Denote the tangles adjacent to ${\cal T}_0$ at $\gamma^*_i$ by ${\mathscr T}_i$  with $i=\{1, \dots, \mu \}$.  Since each essential circle $\gamma^*_i$ is $\phi$-generic, its $\phi$-orbit has $q$ images and the $\phi$-orbit of each tangle ${\mathscr T}_i$  has also $q$ images in the set  $\{{\mathscr T}_i : i=1, \dots, \mu \}$. Using flypes as long if necessary, we modify 
 these tangles so that the orbit of the tangle ${\mathscr T}_i$ is the set ${\mathscr T}_i, \phi ({\mathscr T}_i),  \dots \phi ^{q-1}({\mathscr T}_i)$. Therefore, the union of these modified diagrams constitutes an invariant subset under $\phi$. This means that for a $q$-action whether free or not, the $\mu$ tangles adjacent to ${\cal T}_0$ exhibit a $q$-symmetry induced by  the $q$-rotation $\phi$ on $\mathbb S$. \\

We now consider a tangle $T$ formed by $ \mu \over q$ consecutive tangles of ${\mathscr T}_i$ regrouped by the TBD $\mathcal T_0$ such that the $q$ tangles $ \phi(T) \cup  \phi^2(T) \dots  \cup \phi^q(T)$ contain the tangles adjacent to $\mathcal T_0$. The tangle $T$ will be called the {\bf primary tangle} of ${\cal T}_0$.

Note that each crossing of $\Pi$ is either a crossing in one of the $\mu$ tangles ${\mathscr T}_i, i=1, \dots, \mu$ or is a visible crossing of the main TBD ${\cal T}_0$.\\

However, the behaviour of the visible crossings is not the same in the case of free and non-free $q$-action. This is described in {\it Case 2(a)} for the non-free $q$-periodicity and in {\it Case 2(b)} for the free $q$-periodicity.
\begin{enumerate}
\item {\it Case 2(a)}. {\bf $\Phi$ is a $q-$periodic map of $(S^3,K)$}. It is the Visibility Theorem 5.3 in \cite{co2}.
 \item {\it Case 2(b)}. {\bf $\Phi$ is a free $q$-periodic map of $(S^3,K)$}.\\
By identifying a freely $q$-periodic knot $K$ with a $(q,r)$-lens knot, the map $S_r$ of $\Phi_{q,r}= S_r \circ R_q$ on $S^3$ is nontrivial, which is revealed by the fact that the weight $w$ is not a multiple of $q$. As described in ($\blackdiamond \blackdiamond$), whatever the considered alternating projection $\Pi$ of $K$, the flypes located on the main TBD ${\cal T}_0$ cannot be reduced to the identity map.\\
Using flypes if necessary, the best situation we can obtain is an alternating projection $\Pi$ where $\Phi_\Pi$ is expressed by Equation (19):
$$ \Phi_\Pi= \phi \circ F_0$$ 
where \\
-$\phi$ is conjugate to a rotation $\phi^*$ of order $q$ around an axis $l_2$ such that the simple closed curve $l_1$, core of the main TBD and $l_2$ are linked together exactly once, and \\
-$F_0$ is the composition of standard flypes all occurring in the main TBD ${\cal T}_0$ of $\Pi$.
\end{enumerate}
\end{enumerate}

{\bf Summary.} The analysis in {\it Case 2(b)} constitutes the proof of Theorem 5.4.\\

The Visibility $q$-periodicity Theorem 5.3 implies the following result:

\begin{theorem} (Periodic Projection Theorem) Let $K$ be a prime alternating $q$-periodic knot for $q \geq 3$. Then
in the case $Fix (\widetilde \varPhi)$ is a vertex representing a TBD, there is a {\bf 2}-tangle $T$ and a 2-braid tangle $\sigma^{k}$ such that $K$ has an alternating projection $\Pi$ of the form: 
\begin{equation}
\Pi=\widehat { \sigma^{k} T\,\,   \sigma^{k} T \dots  \sigma^{k} T}.
\end{equation}
where $\sigma$ is a generator of $B_2$.
\end{theorem}

For the case of a free $q$-periodicity, the proof of the Visibility Theorem 5.4 is done  by focusing on the main TBD of any alternating projection.
This gives rise to a projection whose $q$-symmetry is as visible as possible. The $q$-symmetry in the set of the tangles adjacent to ${\cal T}_0$ is realized but it is not possible to realize the $q$-symmetry on $\mathbb S$ for the visible crossings of ${\cal T}_0$. \\ 
The equivalent of Theorem 5.5 for freely periodicity is the following theorem:
 \begin{theorem}(Freely Periodic Projection Theorem)  Let $K$ be a prime alternating freely $q$-periodic knot. Then there is a
 {\bf 2}-tangle $T$ and a 2-braid tangle $\sigma$ such that there is an alternating projection $\Pi$ of $K$  which is the closure of
  \begin{equation}
  T^q \sigma^w
  \end{equation}
 \end{theorem}
 
The proof of Freely Periodic ProjectionTheorem 5.6 is given in the following subsection.

 \subsection{ $q$-periodic and free $q$-periodic alternating projections.}
 
Let $K$ be an alternating knot endowed with a $q$-action. Thanks to the description of a $q$-action on the essential decomposition of $(S^2, \Pi)$, we thus focus on the case where $Fix (\widetilde \varPhi)$ is a vertex associated to a TBD. 
We now introduce the notion of virtual visibility which concerns both the cases of free and non-free actions. 

\begin{definition} Let $K$ be an alternating knot in $S^3$. An alternating projection $\Pi$ of $K$ is {\bf virtually $q$-visible} if there is a tangle $T$, the primary tangle of $\Pi$, and a 2-braid $\sigma$ such that $\Pi$ is the closure 
 
 \begin{equation}
\Pi =\widehat { \sigma^{i_1} T\,\, \sigma^{i_2}  T \,\, \sigma^{i_3} \dots \sigma^{i_q} T}
 \end{equation}
where the occurrences of $T$ and $\sigma$ (with its exponent) alternate in the expression (22).
 \end{definition} 
 
\begin{figure}[h!]
\centering
\includegraphics[scale=.6]{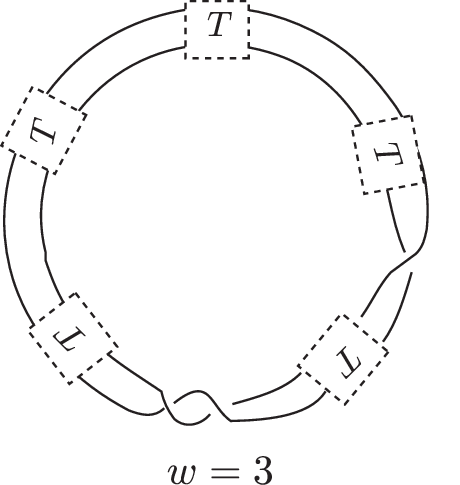}
\caption{A 5-virtually visible projection}
\end{figure}

 \begin{example} Fig.12 exhibits a 5-virtually visible projection which is the closure of $\sigma^2 \,  T \, \sigma \,  T^4$. The main TBD is of weight 3 and valency 5.
\end{example}
 
\begin{figure}
\centering
\includegraphics[scale=.25]{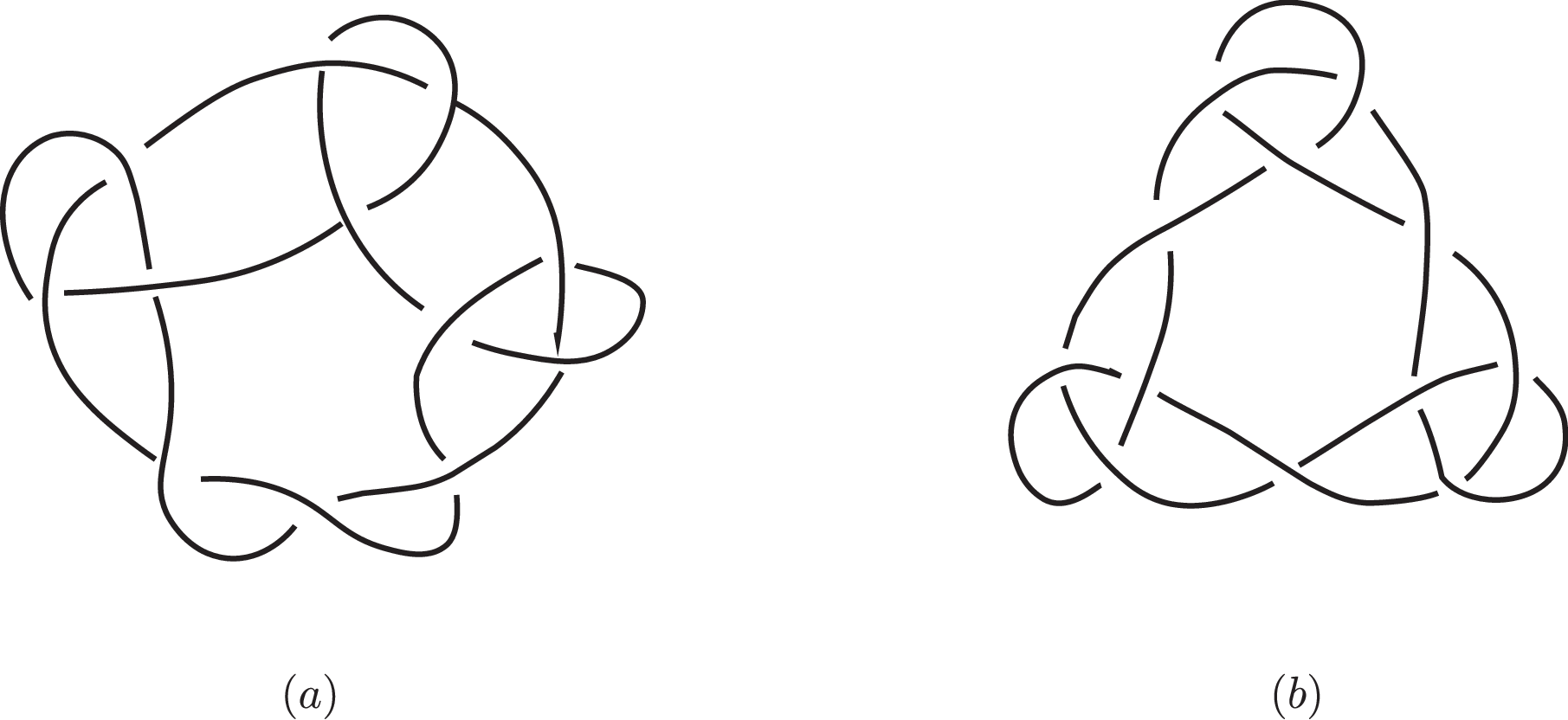}
\caption{Knot $12_{503}$ with {\it (a)} a virtual 3-visibility and {\it (b)} a 3-periodic projection}
\end{figure} 

\begin{example}
 Fig.13 depicts the knot $12_{503}$ in a virtual 3-visible projection $(a)$ and in a 3-periodic projection $(b)$.
\end{example}

In the case of an {\bf alternating} freely periodic knot, we have: 

\begin{proposition} Let $K$ be an alternating freely $q$-periodic knot which is not a torus knot. Then $K$ can be expressed as 
\begin{equation}
K= \widehat{{\rm T}^q \sigma^{2r}}.
\end{equation}
where $\rm T$ is a {\bf 2}-tangle.
\end{proposition}

\begin{notation}  $T$ (italicized T) denotes a primary tangle and $\rm T$ (romanized T) denotes a fundamental lens-tangle  which is described in Proposition 5.8.
\end{notation}
 \begin{figure}[h!]
\centering
\includegraphics[scale=.2]{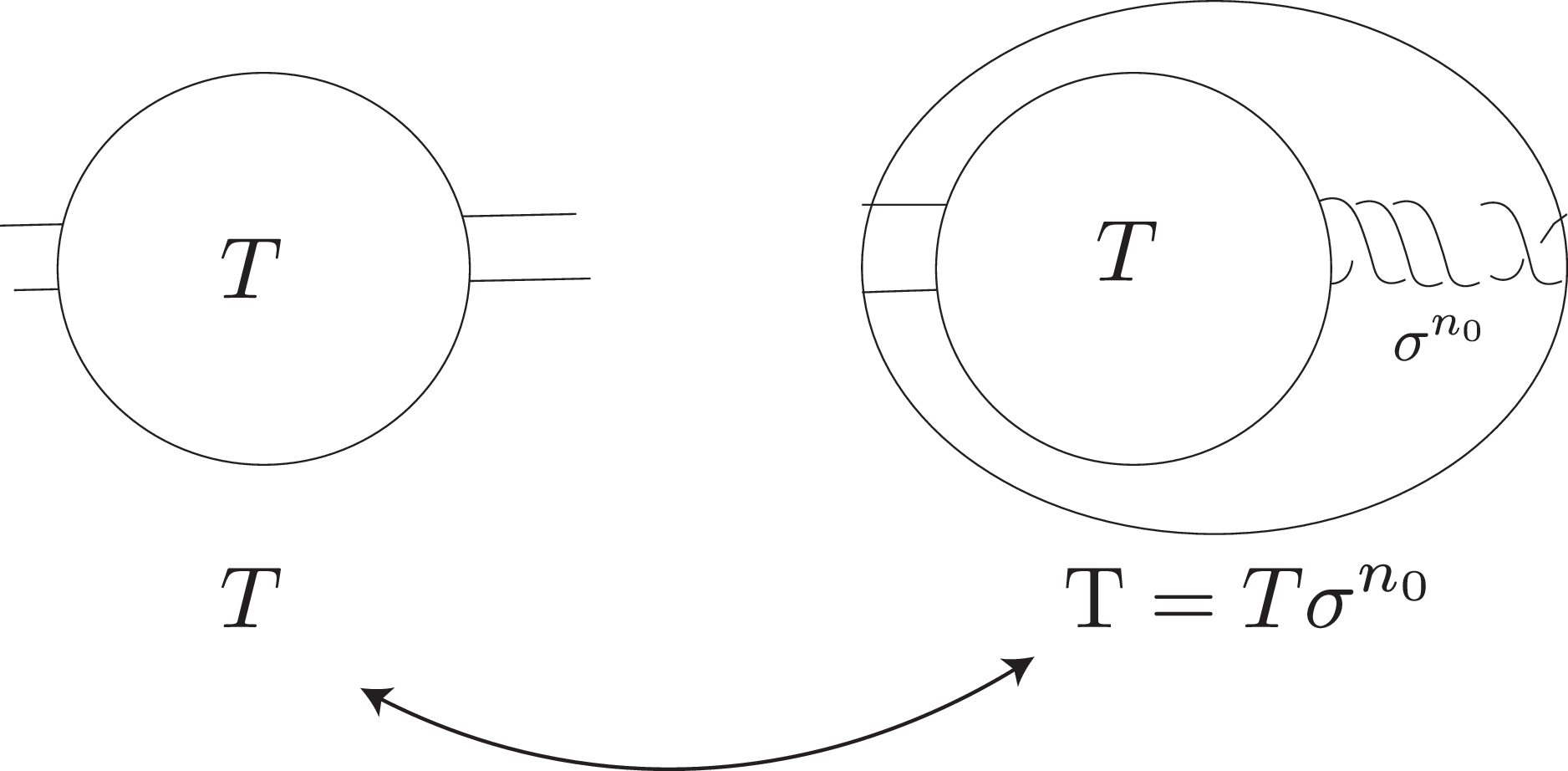}
\caption{Relation between the primary tangle $T$ and the fundamental lens-tangle $\rm T$}\
\end{figure}
\begin{proof} (Proof of Proposition 5.8).
Let $K$ be a a freely $q$-periodic alternating knot which is not a torus knot. Let $\Pi$ be an alternating projection of $K$ with its main TBD ${\cal T}_0$ and its primary tangle $T$.  Since the freely $q$-periodic alternating knots are identified to the $(q,r)$-lens knots, the Visibility theorem 5.4  implies that the 2-version of $\lambda (\rm T)^q \, \sigma^{2r}$ is ${\rm T}^q \, \sigma^{2r}$ and  that $\rm T$ is necessarily a {\bf 2}-tangle obtained from a sum of the (alternating) primary tangle $T$ with a power of $\sigma$. 
\end{proof}

We now give the proof of Freely Periodic ProjectionTheorem.
\begin{proof} (Proof of Theorem 5.6).
By Proposition 5.8 the freely $q$-periodic knot $K$ can be represented as a $(q,r)$-lens knot with the expression of Equation (23).
If the closure of ${\rm T}^q \, \sigma^{2r}$ is alternating, we have $T={\rm T}$ and $w=2r$.\\
Otherwise, as described in the proof of Proposition 5.8, we have ${\rm T}=T \sigma^{-n_0}$ for some $n_0$ and we can express $K$ as  
\begin{equation}
K= \widehat {(T \sigma^{-n_0}T \sigma^{-n_0} \dots T \sigma^{-n_0}) \sigma^{2r}}
\end{equation}
(Fig.14 describes how the fundamental lens-tangle $\rm T$ is related to the primary tangle $T$).

Performing flypes and using Reidemeister moves of type II if necessary, we obtain an alternating projection $\Pi$ given by Equation (21)
$$\Pi = {\widehat {T^q {\sigma}^w}}.$$
The main TBD ${\cal T}_0$ is of valency $q$ and weight $w$ with $|w| \leq |2r|$. 
\end{proof}

Let $\Pi$ be ${\widehat {T^q {\sigma}^w}}$ where $T$ is its primary tangle and $w$ is an odd integer which is not a multiple of $q$.

How to find the fundamental lens-tangle from the primary tangle and how to find the coefficient  $r$ of the pair $(q,r)$ of free visibility of $K$?

According to the parity of $w$, we have two cases: 
\begin{enumerate}
\item if $w=2r$, we are done. Then $K$ is a $(q, {w \over 2})$-lens knot. The fundamental lens-tangle $\rm T$ is the primary tangle $T$.

\item if $w$ is odd, the coefficient $r$ of the $(q,r)$-visibility of $K$ is obtained by solving the following diophantine equation with unknowns $n_0$ and $r$:
\begin{equation}
q {\bold n_0}+w=2{\bold r}
\end{equation}
The assumption that $w$ is odd implies that $q$ and  $n_0$ are both odd. Each solution $(n_0,r)$ of Equation (21) produces a $(q,r)$-visibility for $K$. The strategy is to extend the primary tangle $T$ to the tangle ${\rm T}=T \sigma^{n_0}$ such that $n_0$ is of the same sign than $w$. Then $K$ is the closure ${\widetilde \Pi}=\widehat { {\rm T}^q \sigma^{2r}}$ where $\rm T$ is the alternating fundamental lens-tangle of its $(q,r)$-visibility. 
\end{enumerate}

 In what follows, we will establish some restrictions on the integers that are possible periods of freely periodic alternating {\bf knots}.

\begin{figure}[h!]    
   \centering
    \includegraphics[scale=.35]{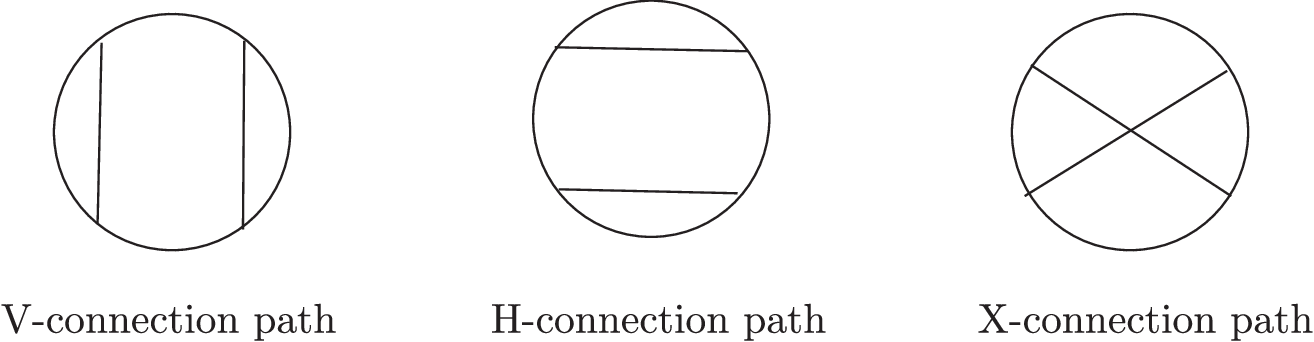} 
\caption{Connection paths of a tangle}       
\end{figure}

\begin{proposition} Let $K$ be a prime oriented alternating knot in $S^3$ with a free $q$-action. Then $q \geq 3$ is an odd integer.
\end{proposition}

We first recall that for each tangle $T$, the projection $\Pi \cap T$ connects the 4 endpoints  NW, NE, SE and SW in pairs with 3 possible {\bf connection paths} (Fig.15).\\
If $\Pi$ connects\\
(1) NW to NE and SW to SE, $T$ is a H-connection path.\\
(2) NW to SW and NE to SE, $T$ is a V-connection path.\\
(3) NW to SE and SW to NE, $T$ is a X-connection path.

\begin{proof} Equation (20) means that $\rm T$  is a ${\bf 2}$-tangle and that $\omega^{2r}$ corresponds to the $2r$ visible crossings of the band. Therefore, for $K$ to be a knot, $\rm T$ must be a X-connection path and $q \geq 3$ must be an odd integer.                                                                                                                                                                         
\end{proof}

We have the following proposition which describes some cases of alternating projections from which one can easily deduce that they admit or not a $q$-action.

Assume that $q \geq 3$ is an odd integer.\\
\begin{proposition}
Let $K$ be a prime oriented alternating knot in $S^3$ and $\Pi$ be an alternating $q$-virtually visible projection of $K$ such that 
$K$ is the closure of the tangle in Equation (22):
$$\sigma^{i_1} \, T\, \sigma^{i_2} \dots \sigma^{i_q}\, T$$
where the primary tangle $T$ is symmetric.\\
The main TBD has its valency equal to $q$ and its {\bf visible weight} equal to $w=i_1+\dots +i_q$.
\vspace{-.3cm}
 \begin{enumerate}
\item  if $w$ is a multiple of $q$, $K$ is $q$-periodic. 
 \item if $w$ satisfies $gcd(w,q)=1$, $K$ is freely $q$-periodic. 
  \end{enumerate}
\end{proposition}

\begin{proof}
$(a)$ $w= a \, q$. Since $\rm T$ is symmetric i.e.  $\rm T \sim_f \rm F(T)$ (Definition 3.10), by making flypes if necessary, we can obtain the 
configuration $\sigma^a \, T \, \sigma^a \dots \sigma^a \, T$.

$(b)$ By performing flypes, we can obtain the tangle $T^q \sigma^{w}$ as described in \S 5.5. The proof of Proposition 5.8 explains how to reveal the $(q,r)$-visibility of $K$.
\end{proof}

\begin{example}
(a) The knot $10_{75}$ (Fig.16) is a $(3,1)$-lens knot. 

\begin{figure}[h!]
\centering
\includegraphics[scale=.6]{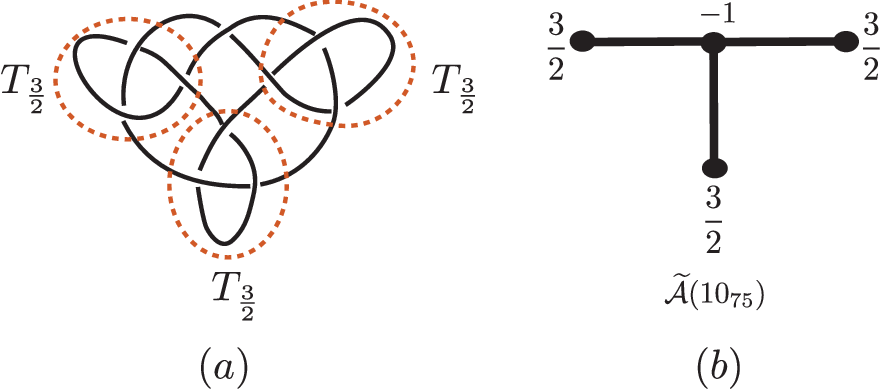}
\caption{Knot $10_{75}$ with its $(a)$ virtual 3-visibility and its $(b)$ essential structure tree}
\end{figure} 

The essential structure tree ${\widetilde {\cal A}} (10_{75})$ has four vertices such that one corresponds to a TBD ${\mathcal T}_0$ of weight $w=-1$ and each of the other three corresponds to the tangle $T_{3 \over 2}$.

\begin{figure}[h!]
\centering
\includegraphics[scale=.7]{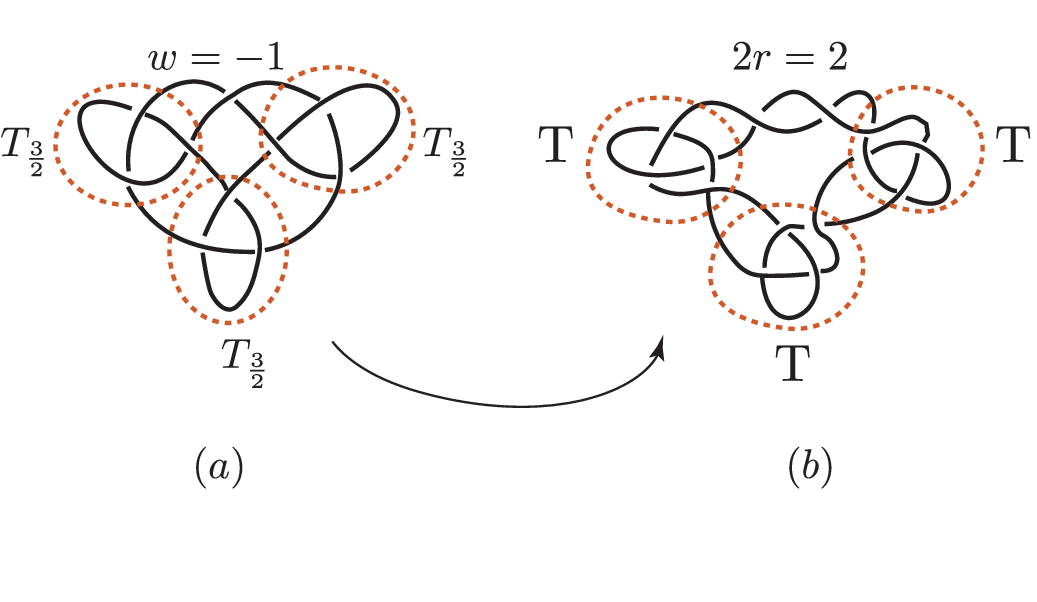}
\caption {How to obtain a $(3,1)$-lens visible projection from a 3-virtually visible projection of $10_{75}$}
\end{figure} 
 The knot $10_{75}$ has a 3-virtually visible projection with a TBD of visible weight $w=-1$.
In this case, Equation (17)  is therefore equal to:
\begin{equation}
3 {\bold n} -2{\bold r}=1
\end{equation}
Since $gcd(2,3)=1$, we have at least a solution.
The pair $(n,r)=(1,1)$ gives rise to the form:
\begin{equation}
10_{75}=\widehat{\rm T^3 \sigma^2}.
\end{equation}

\begin{figure}[h!]  
\centering
\includegraphics[scale=.2]{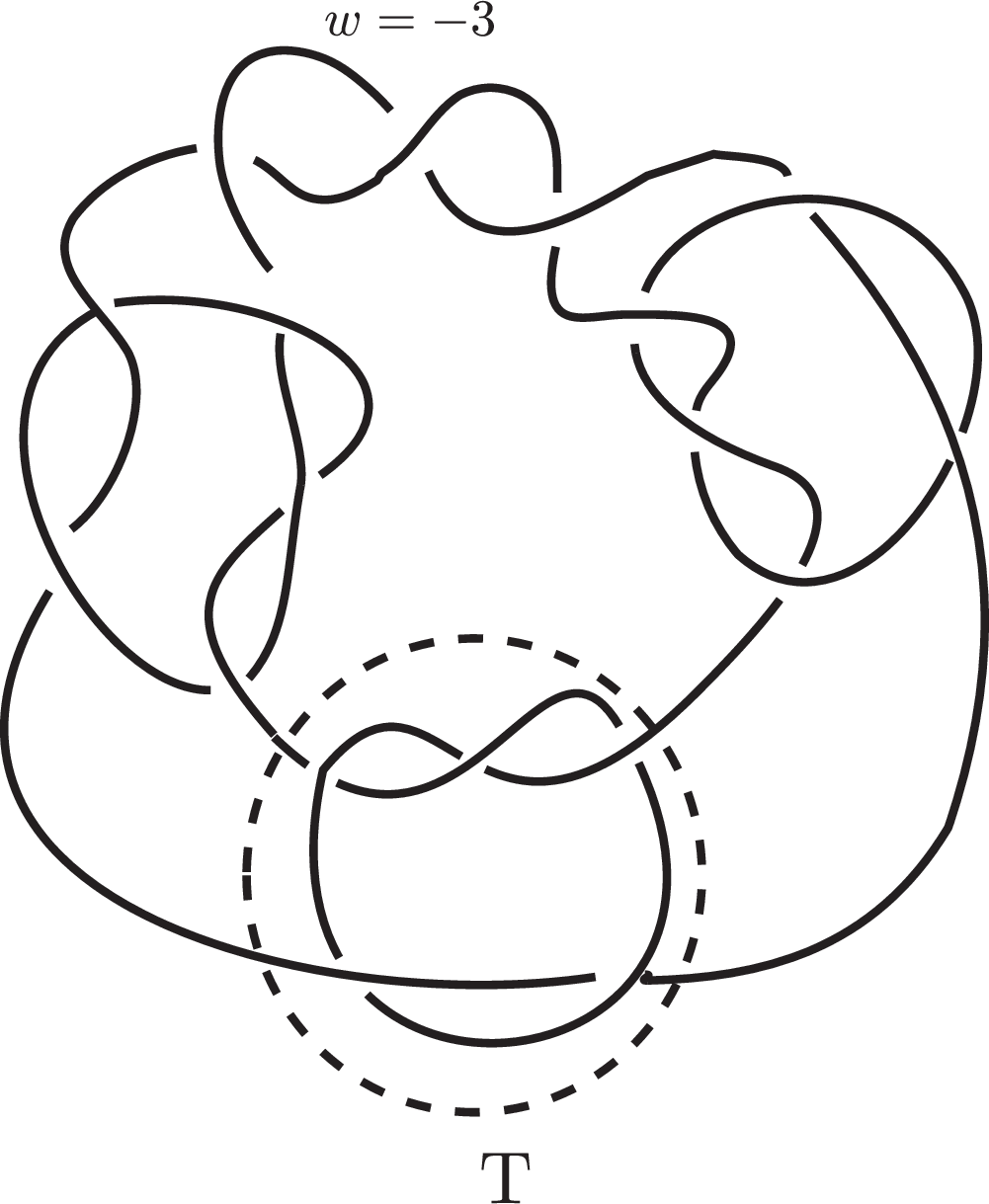}
\caption{A virtually 3-periodic 18-crossing knot}
\end{figure} 

(b) Consider the 18-crossing alternating knot depicted in Fig.18. The described projection is virtually 3-periodic with visible weight $w=-3$. However, since the tangle $\rm T$ is not symmetric, we cannot deduce any 3-periodic projection. By \cite{co2}, we conclude that the knot  is not 3-periodic.
 \begin{figure}[h!]  
\centering
\includegraphics[scale=.2]{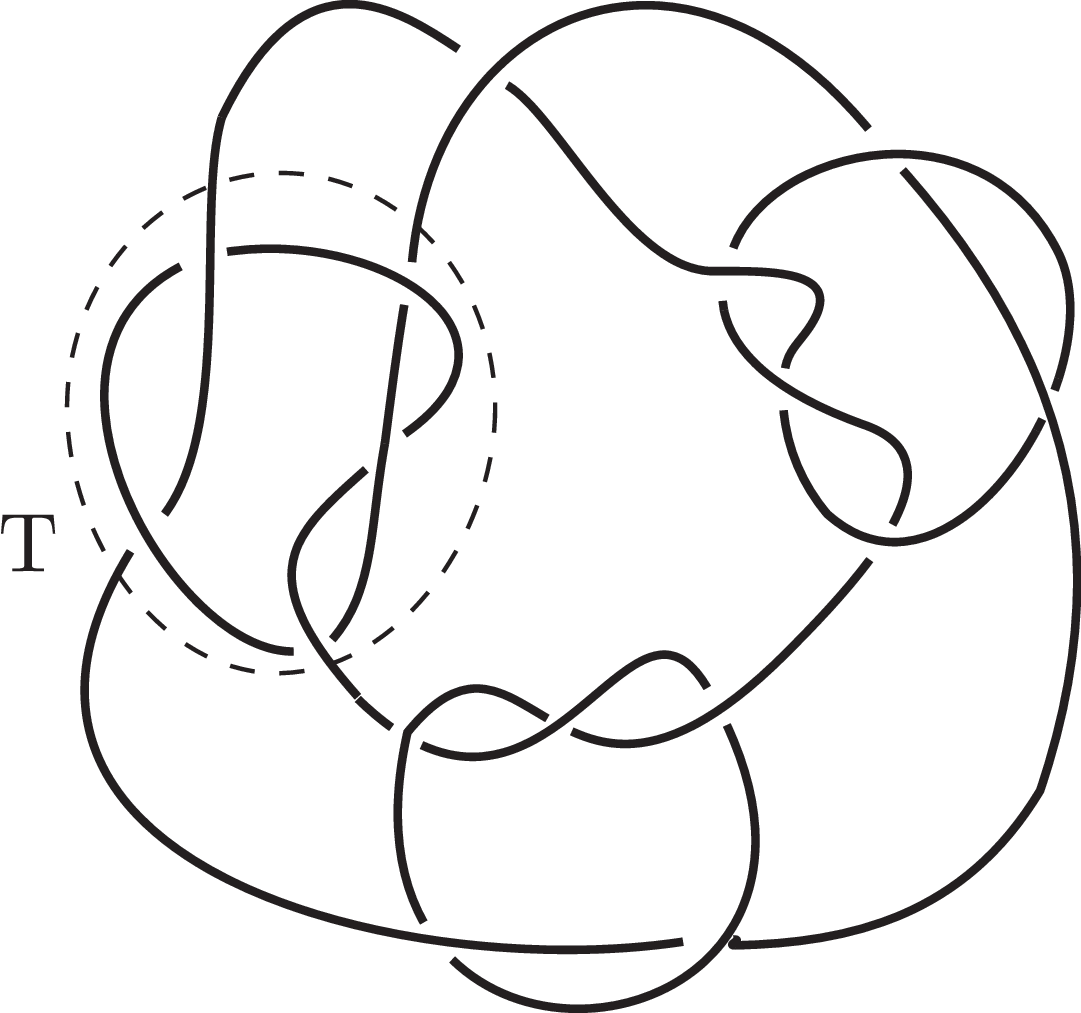}
\caption{A virtually 3-periodic 16-crossing knot}
\end{figure} 

(c) Consider the 16-crossing alternating knot described in Fig.19.
This knot is not 3-periodic since the visible weight is not a multiple of 3. Since the tangle $\rm T$ is not symmetric,  it is not freely 3-periodic either. 
 \end{example}
 
 \section{Murasugi decomposition and $q$-action on alternating knots}
From \cite{co2}, to any $q$-periodic prime oriented alternating knot $K$, we can exhibit a $q$-periodic alternating projection.  This implies that the adjacency graph ${\cal G}(K)$ of the Murasugi decomposition into atoms of $K$ admits an automorphism of order $q$ (\cite{co}, see Theorem 6.1 below). In this section, we will show that the Murasugi decomposition on a freely $q$-periodic alternating knot enables to extract  some interesting informations on its free $q$-periodicity. The strategy is to locate the position of the Murasugi atoms with respect to the essential decomposition by the essential Conway circles into jewels and TBDs of a freely $q$-periodic alternating knot.
 
 \begin{figure}[h!]  
\centering
\includegraphics[scale=.3]{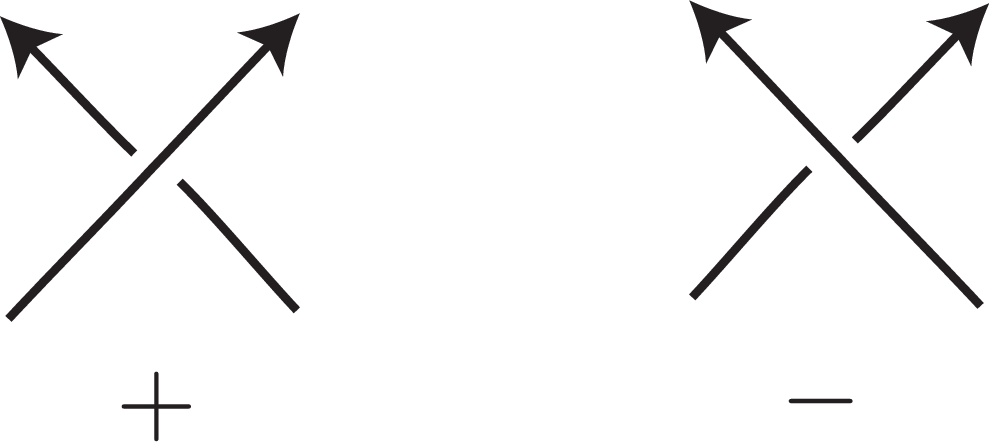} 
\caption{Sign of a crossing}
\end{figure} 

Let $L$ be a prime {\bf oriented} alternating link. Consider an alternating projection $\Pi$ of $L$ endowed with the orientation of $L$. By splicing  $\Pi$ at each crossing while respecting the orientation, we get a bunch of {\it Seifert circles} which are oriented, simple and disjoint circles. These circles bound {\it Seifert disks} on the projection sphere. By connecting these disjoint disks by half-twists that correspond to crossings, we obtain an {\it algorithmic Seifert surface} $\widetilde \Sigma$ based on $\Pi$ whose boundary is $L$. We can express $\widetilde \Sigma$ as the {\it diagrammatical Murasugi sum} of the Seifert surfaces of the atoms of $L$. As done in (\cite{quwe1}, \cite{stoi}), an atom is a prime oriented special alternating link (where {\bf special} means that all of its crossings are of the same sign). So we have positive and negative atoms. The collection of $L$-atoms ${\cal C}_a (L)$ is independent of the choice of a projection of $L$, provided it is alternating and is therefore an isotopy invariant.\\
We define the {\bf adjacency graph} ${\cal G}(L)$ of an oriented prime alternating $L$ as a bipartite tree whose vertices represent negative or positive atoms; two vertices of ${\cal G}(L)$ are connected if and only if their two associated atoms with opposite sign share a common Seifert circle. The adjacency graph is an isotopy invariant in the class of oriented alternating links (\cite{quwe1}).

Let $K$ be a prime oriented alternating knot endowed with a $q$-action generated by a $q$-map $\Phi$ (where $q \geq 3$ is an odd integer). The map $\Phi$ induces an automorphism $\widetilde {\varPhi}$ on the essential structure tree $\widetilde {\cal A} (K)$ and an automorphism on the adjacency graph $\cal{G}(K)$. 

In the previous papers (\cite{co}\cite{co2}), for a non-free $q$-action, we showed:

\begin{theorem}
Let $K$ be a prime alternating knot endowed with a non-free $q$-periodic map $\Phi$ ($q \geq 3$ is an odd integer). Consider its collection of atoms ${\cal C}_a(K) = \{L_1, L_2, \dots,L_s \}$ and its adjacency graph $\mathcal{G}(K)$. Then the labelled graph $\cal{G}(K)$ has an automorphism ${\mathscr F}$ of order $q$ induced by $\Phi$. Either each atom $L_i$ of $K$ is $q$-periodic or it occurs a multiple of $q$ times in ${\cal C}_a(K)$.
\end{theorem}
{\bf Question}: What is the equivalent of Theorem 6.1 for a free $q$-action?

{\bf Answer} Indeed with the same arguments, we have a similar result for the case of the non-free periodicity.

We first recall that if $Fix (\widetilde {\varPhi})$ corresponds to a jewel $J_0$, the knot $K$ is a fortiori $q$-periodic (Proposition 5.5). Since there are no flypes in $J_0$, the $q$-periodicity is visible and any constituent atom of $J_0$ is $q$-periodic. 

 For a free $q$-action, we have the following result:
\begin{theorem}
 Let $K$ be a prime alternating knot endowed with a free $q$-action given by $\Phi$ (where $q \geq 3$ is odd). Consider its collection of atoms ${\cal C}_a(K)$ and its adjacency graph $\mathcal{G}(K)$. Then\\
 1) the labelled graph $\mathcal{G}(K)$ has an automorphism ${\mathscr F}$ of order $q$ induced by $\Phi$ such that $Fix \,{\mathscr F}$ corresponds to a unique $q$-freely periodic atom and \\
 2) each of the other atoms of ${\cal C}_a(K)$ appears a multiple of $q$ times in ${\cal C}_a(K)$.
\end{theorem}
For the proof of Theorem 6.2, we need the following lemma which uses the essential decomposition of $(S^2, \Pi)$ into jewels and TBDs.
\begin{figure}[h!]  
\centering
\includegraphics[scale=.7]{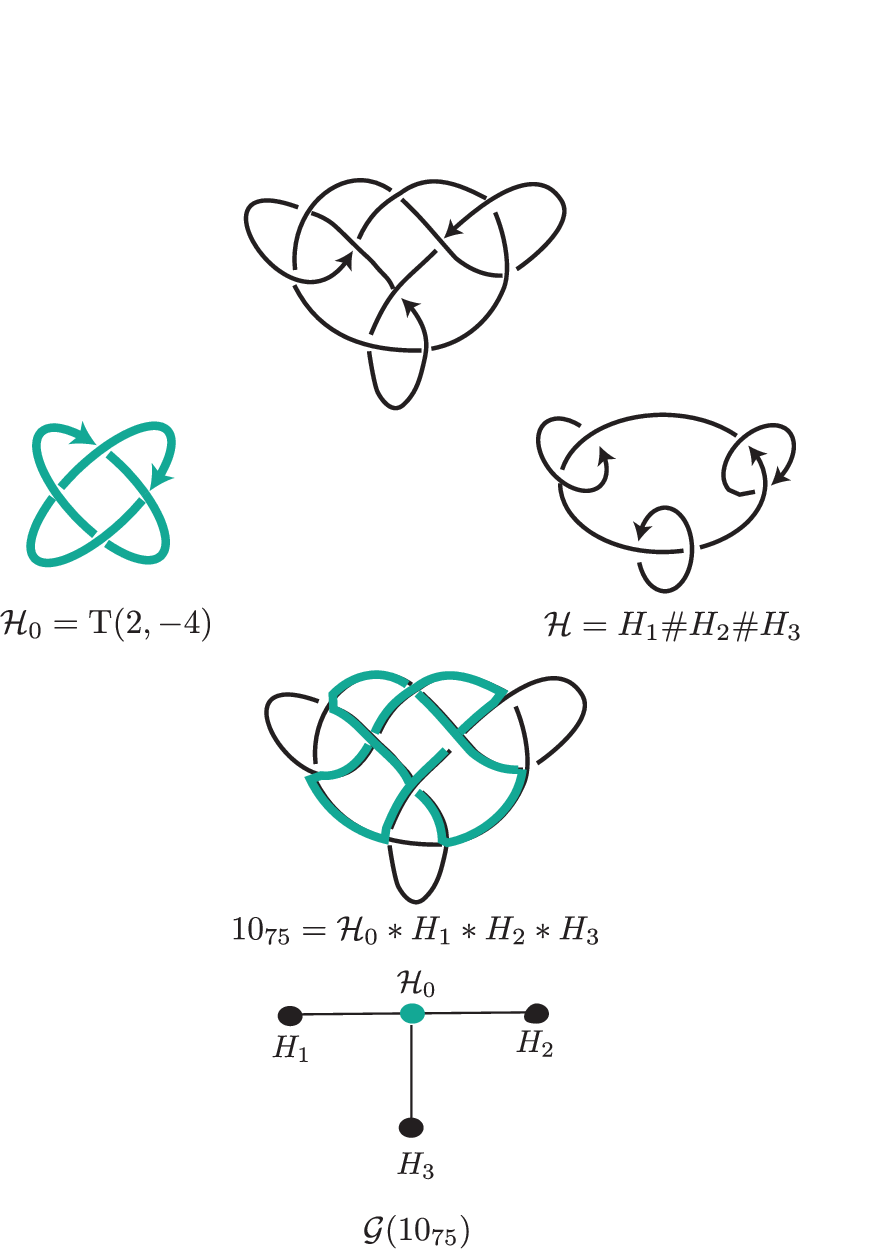}
\caption{Murasugi decomposition of $10_{75}$ and its adjacency graph}
\end{figure} 

\begin{lemma} Let $q \geq 3$ be an odd integer. Let $K \subset S^3$ be a prime alternating knot endowed with a free $q$-action generated by $\Phi$. Consider the collection of $K$-atoms ${\cal C}_a(K) $ and its adjacency graph $\mathcal{G}(K)$. Let ${\mathscr F}$ be the automorphism of order $q$ induced by $\Phi$ on $\mathcal{G}(K)$. Then
$Fix\,{\mathscr F}$ is a vertex corresponding to a unique atom.
\end{lemma}
We call the atom corresponding to the vertex $V_0=Fix {\mathscr F}$, the {\bf main atom} of $K$.

\begin{proof} (proof of Lemma 6.2)
Let $\Pi$ an alternating projection of $K$.
Since $K$ is freely $q$-periodic, by Proposition 5.7, $Fix \,(\widetilde {\varPhi})$ is a vertex corresponding to a TBD: the main TBD ${\cal T}_0=(\Sigma, \Gamma=\Pi \cap \Sigma)$ of $\Pi$. \\
The Seifert algorithm makes a set of Seifert circles appear from an alternating projection $\Pi$. 
Looking at how the arcs of $\Gamma$ are related to this bunch of Seifert circles, we deduce that the 1-manifold $\Gamma=\Pi \cap \Sigma$ of the TBD ${\cal T}_0$ belongs to a same atom of ${\mathcal C}_a(K)$; this corresponds to  a vertex $V_0$ of $\mathcal{G}(K)$ (see \cite{quwe1}). Therefore, we have: $Fix \,(\widetilde {\varPhi})= {\cal T}_0$ implies that $Fix \,{\mathscr F}=V_0$.\\
If $Fix \,{\mathscr F}$ were not reduced to a unique vertex, $Fix \,(\widetilde {\varPhi})$ would correspond to a jewel, which implies that the involved knot would be $q$-periodic (Proposition 5.6). \\
\end{proof}

\begin{example}  Fig.21 depicts the Murasugi decomposition of the knot $10_{75}$ where  $H_1, H_2, H_3$ are three  positive Hopf links and ${\cal H}_0$ is a negative atom which is the torus link $ \rm T (2,-4)$. The central vertex of the adjacency graph  ${\cal G } (10_{75})$  represents ${\cal H}_0$. The 1-manifold $\Gamma$ of the main TBD ${\cal T}_0=(\Sigma, \Gamma=\Pi \cap \Sigma)$ belongs to ${\cal H}_0$ the main atom of $10_{75}$.
\end{example}

\begin{proof} (Proof of Theorem 6.2)
The automorphism ${\mathscr F}$ on the adjacency tree $\cal{G}(K)$ is of order $q$ with $Fix \, {\mathscr F} = V_0$. Each vertex distinct of $V_0$ has its generic orbit under ${\mathscr F}$ and $V_0= Fix \,{\mathscr F}$ corresponds to the free $q$-periodic atom responsible of the free $q$-periodicity of $K$.
\end{proof}

Thus, we have the following nice corollary where $q$ is an odd integer $>2$.
\begin{corollary} Let $K \subset S^3$ be a prime alternating $q$-freely periodic knot with its main atom $K_0$. Then $K$ and its main atom $K_0$ have the same coefficient $r$ for their $(q,r)$-visibiliy.
\end{corollary}

In the following example, we show how to determine the coefficient $r$ by focusing  on the main atom of $K$.

\begin{example} We first apply the Murasugi decomposition on $10_{75}$. 
 Since ${\cal H}_0={\rm T}(2,-4)$ is the main atom of $10_{75}$, it is sufficient to derive the coefficient $r$ from the torus link ${\rm T}(2,-4)$.\\
We now apply Lemma 4.2 to ${\rm T}(2,-4)$. The number $rm-n= -4r-2$ is divisible by $q=3$ for $r=1$. Thus, $10_{75}$ is a $(3,1)$-lens knot by Corollary 6.1.
\end{example} 

\section{Addendum}
There is some overlap between the results of K. Boyle's paper ( \cite{boy}) and this one, as well as our previous paper (\cite{co2}).The Flyping theorem is central but the techniques differ with somewhat different visualizations. For example, the symmetry of freely $q$-periodic alternating knots in our approach makes it easy to identify them as $(q,r)$-lens knots.

\end{document}